\documentclass[twoside, 12pt]{amsart}
\usepackage[active]{srcltx} 
\usepackage[all]{xy}
\CompileMatrices
\usepackage{ifthen}
\usepackage{amsfonts}
\usepackage{amssymb}
\usepackage{amsthm}
\usepackage{amsmath}
\usepackage{url}
 \newtheorem{thm}{Theorem}[section]
 \newtheorem{prop}[thm]{Proposition}
 \newtheorem{cor}[thm]{Corollary}
 \newtheorem{lem}[thm]{Lemma}
  \newtheorem{ques}[thm]{Question}
 
 \newtheorem*{bthm}{Theorem}

\theoremstyle{definition}

\newtheorem{ex}[thm]{Example}

\theoremstyle{remark}
\newtheorem{rem}[thm]{Remark}

  1

\newcommand{\fp}{\ifmmode {\mathbb{F}_p}\else$\mathbb{F}_p$\ \fi}
\newcommand{\zp}{\ifmmode {\mathbb{Z}_p}\else$\mathbb{Z}_p$\ \fi}

\newcommand{\zpMod}{\ifmmode\mbox{$\zp$-Mod}\else$\zp$-Mod \fi}
\newcommand{\Mod}{\ifmmode\mbox{$\Lambda$-Mod}\else$\Lambda$-Mod \fi}
\renewcommand{\mod}{\ifmmode\mbox{$\Lambda$-mod}\else$\Lambda$-mod
\fi}
\newcommand{\La}{\ifmmode\Lambda\else$\Lambda$\fi}

\newcommand{\Hom}{{\mathrm{Hom}}}
\newcommand{\Ext}{{\mathrm{Ext}}}

\newcommand{\tor}{{\mathrm{tor}}}
\newcommand{\E}{{\mathrm{E}}}

\newcommand{\rk}{{\mathrm{rk}}}

\renewcommand{\H}{\mathrm{H}}

\newcommand{\M}{\ifmmode {\frak M}\else${\frak M}$ \fi}
\newcommand{\m}{\ifmmode {\frak m}\else$\frak m$ \fi}
\newcommand{\p}{\ifmmode {\frak p}\else${\frak p}$\ \fi}
\renewcommand{\P}{\ifmmode {\frak P}\else${\frak P}$\ \fi}
\newcommand{\e}{\ifmmode {\mathcal{E}}\else$\mathcal{E}$ \fi}
\newcommand{\G}{\ifmmode {\mathcal{G}}\else${\mathcal{G}}$\ \fi}
\renewcommand{\d}{\ifmmode {\mathcal{ D}}\else${\mathcal{D}}$\ \fi}
\newcommand{\A}{\ifmmode {\mathcal{A}}\else${\mathcal{ A}}$\ \fi}

\renewcommand{\projlim}[1] {{\lim\limits_{\stackrel{\displaystyle
\longleftarrow}{#1}}}}

\newcommand{\dirlim}[1]
{{\lim\limits_{\stackrel{\displaystyle \longrightarrow}{#1}}}}
\renewcommand{\in}{\ \epsilon\ }

\newcommand{\kl}{[\![}
\newcommand{\kr}{]\!]}
\newcommand{\Qp}{\ifmmode {{\Bbb Q}_p}\else${\Bbb Q}_p$\ \fi}
\newcommand{\qp}{\ifmmode {{\Bbb Q}_p}\else${\Bbb Q}_p$\ \fi}
\newcommand{\Q}{\ifmmode {\Bbb Q}\else${\Bbb Q}$\ \fi}

\newcommand{\Ind}{\mathrm{Ind}}

\newcommand{\Ref}{\mbox{$\mathcal{R}$}}

\newcommand{\W}{\mbox{$\mathcal{W}$}}

\begin{document}

\title[A noncommutative Weierstrass preparation theorem]{A noncommutative Weierstrass preparation theorem and applications to Iwasawa theory}%
\author{Otmar Venjakob}%
\address{Universit\"{a}t Heidelberg\\ Mathematisches Institut\\
Im Neuenheimer Feld 288\\ 69120 Heidelberg, Germany.}
\email{otmar@mathi.uni-heidelberg.de}
\urladdr{http://www.mathi.uni-heidelberg.de/\textasciitilde
otmar/}
\thanks{During this research, the author
 has been supported by the EU Research Training Network ``Arithmetical Algebraic Geometry"
 at the Department of Pure Mathematics and Mathematical Statistics, Cambridge.}



\date{11.03.02}%
\begin{abstract}
 In this paper and a forthcoming joint one  with Y. Hachimori \cite{hachi-ven} we study Iwasawa
 modules over an infinite Galois extension $k_\infty$ of a number field $k$ whose Galois group
 $G=G(k_\infty/k)$
  is isomorphic to the semidirect product of two copies of the $p$-adic integers \zp.
 After first analyzing some general algebraic properties of the corresponding Iwasawa algebra, we
 apply these results to the Galois group $X_{nr}$  of the $p$-Hilbert class field over
 $k_\infty.$ The main tool we use is a version of the Weierstrass preparation theorem, which we prove
  for certain skew power  series with coefficients in a not necessarily commutative local ring. One striking result in our work is the discovery of the abundance of faithful torsion $\La(G)$-modules,
 i.e.\ non-trivial torsion modules whose global annihilator ideal is zero. Finally we show
 that the completed group algebra $ \fp \kl G \kr  $ with coefficients in the finite field
 $\mathbb{F}_p$  is a unique factorization domain in the sense of Chatters \cite{cha}.
\end{abstract}

\maketitle

\section{Introduction}

One of the highlights  of noncommutative Iwasawa theory so far is
the structure theory of Iwasawa modules due to Coates, Schneider
and Sujatha. In \cite{co-sch-su} they proved that every torsion
$\La(G)$-module decomposes - up to pseudo-isomorphism - into the
direct sum of cyclic modules, where $\La(G)=\zp \kl G \kr$ denotes
the Iwasawa algebra of a pro-$p$-group $G$ belonging to a certain
class of $p$-adic groups, viz the class of $p$-valued ones. In
this paper we want to illustrate some features of their theory in
the possibly easiest noncommutative example: we take as group the
semi-direct product $G=H\rtimes \Gamma$ of two copies of the
$p$-adic integers $\zp\cong H\cong\Gamma$ with a non-trivial
action of $\Gamma$ on $H.$

If $M$ is a finitely generated left module over a ring $R,$ we
define, as usual, $Ann(M)$ to be the set of all $r$ in $R$ such
that $r\cdot M=0.$ If $R$ is commutative and $M$ is a torsion
$R$-module, we always have $Ann(M)\neq0.$ However, when $R$ is
non-commutative, we define a torsion $R$-module $M$ to be faithful
or bounded, according as $Ann(M)=0$ or $Ann(M)\neq0.$ When
$R=\Lambda(G),$ with $G=H\rtimes \Gamma$ as above, both types of
modules occur. But what is very striking is that those modules
which occur most naturally in arithmetic ( see below and
\cite{hachi-ven}), namely the modules which are finitely generated
over $\Lambda(H),$ are either faithful or pseudo-null. In fact, an
even stronger result is true. Let  $\mathcal{C}$ be the
subcategory of all pseudo-null modules in the category
$\La(G)$-mod of all finitely generated $\La(G)$-modules.

\begin{bthm}[Theorem \ref{th1}]

Let $M$ be a $\La(G)$-module which is finitely generated over
$\La(H).$ Then he image of $M$  in the quotient category
$\mod/\mathcal{C}$  is a completely faithful, cyclic object.
\end{bthm}

Here an object $\mathcal{M}$ of the quotient category
$\mod/\mathcal{C}$ is called {\em completely faithful} if all its
nonzero subquotient objects are faithful (see {\S}
\ref{com-faithful}).

In section \ref{tate} we will consider $p$-adic Lie extensions
$k_\infty$ over a number field $k$ whose Galois group
$G=G(k_\infty/k)$ is isomorphic to a semidirect product of the
type studied above. Such examples arise by adjoining the $p$-power
roots of an element $\alpha\in k$ which is not a root of unity to
the cyclotomic \zp-extension of $k,$ assuming that the latter
field contains the $p$th roots of unity $\mu_p,$ or $\mu_4$ if
$p=2.$  If we take for example $\alpha=p$ and $k=\Q(\mu_p),$ we
shall show  that the Galois group $X_{nr}$ of the $p$-Hilbert
class field gives rise to a  completely faithful object.

At this opportunity  we should mention a further result on the
module $X_{nr}.$  The remarkable duality between the inverse and
direct limit  of the $p$-primary ideal class groups in the tower
$k_\infty$ of number fields  which was observed by Iwasawa for the
cyclotomic \zp-extension - compare also with results of McCallum
\cite{mcCal} or Nekovar \cite{nek} in the case of multiple
\zp-extensions -  holds in the following (noncommutative) case

\begin{bthm}[Theorem \ref{xnr}]
Let $k=\mathbb{Q}(\mu_p)$ and
$k_\infty=k(\mu_{p^\infty},p^{p^{-\infty}}).$ Then, there is an
injective homomorphism of $\Lambda(G)$-modules
\[Cl(k_\infty)(p)^\vee\to \E^1(X_{nr})\] with pseudo-null
cokernel, i.e.\ the Pontryagin dual of the direct limit of the
$p$-ideal class groups is pseudo-isomorphic to the Iwasawa-adjoint
of the inverse limit. In particular, if $X_{nr}$ is pseudo-null,
then   the ideal class group $Cl(k_\infty)(p)=0$ vanishes.
\end{bthm}

It is easy to see that $X_{nr}$ is non-zero if and only if $p$ is
an irregular prime. We know no example at present where $X_{nr}$
is not pseudo-null, and it may well be that is always the case.
In a forthcoming paper \cite{hachi-ven} (with Y. Hachimori) we
will prove that certain Selmer groups or rather their Pontryagin
dual are indeed non-zero completely faithful objects in the
quotient category. Thus faithful modules definitely occur in the
world of arithmetic. Hence, in these cases, one cannot hope to
obtain  a link with L-functions
 via a global annihilator or a characteristic ideal, which simply does not exist for such
 completely faithful   modules.

To obtain the algebraic results on modules over $\La(G),$ where
now in general  $G$ is a uniform pro-$p$-group which is isomorphic
to the semidirect product $H\rtimes\Gamma$ of a uniform
pro-$p$-group $H$ with $\Gamma\cong\zp,$ we proceed as follows:
Firstly, we identify $\La(G)$ with the skew power series ring
$=R[[Y;\sigma,\delta]]$ over the subalgebra $R:=\La(H)$ of
$\La(G).$ Secondly, we study in some generality skew power series
rings $A:=R[[Y;\sigma,\delta]]$ the definition of which is given
in section \ref{powerseries} and prove a generalized Division
Lemma and Weierstrass Preparation Theorem  for such rings (see
Theorem \ref{WP} and corollary \ref{cor-WP}). Here we give only
one remarkable  consequence of this preparation theorem.  By
$\mbox{$\La(G)$-mod}^H$ and $\mathcal{C}^H $ we denote the
category of $\La(G)$-modules which are finitely generated as
$\La(H)$-module, and its full subcategory of those pseudo-null
$\La(G)$-modules, which are finitely generated over $\La(H),$
respectively. The skew polynomial ring $R[Y;\sigma;\delta]\cong
R[Z,\sigma],$ $Z=Y+1,$ possesses as a ring of quotients
$Q[Z;\sigma]$ where $Q=Q(H)$ denotes the skewfield of fractions of
$\La(H).$  We prove the following

 \begin{bthm}[Theorem \ref{maxleft-simple}]
The quotient category $\mbox{$\La(G)$-mod}^H/ \mathcal{C}^H $ can
be identified with a full subcategory of the category
$\mbox{$Q[Z;\sigma]$-mod}$ of finitely generated
$Q[Z;\sigma]$-modules \[\mbox{$\La(G)$-mod}^H/ \mathcal{C}^H
\subseteq \mbox{$Q(H)[Z;\sigma]$-mod}.\] Under this
identification,  the equivalence classes of simple objects of the
quotient category correspond uniquely to equivalence classes of
maximal left ideals of $Q[Z;\sigma]$ which contain a distinguished
polynomial. Furthermore, bounded objects correspond to bounded
$Q[Z;\sigma]$-modules.
\end{bthm}

 An analogous statement holds if we replace $Q[Z;\sigma]$ by the skew Laurent polynomial ring
 $Q[Z,Z^{-1};\sigma].$ In the case where $G\cong\zp\rtimes\zp$ with non-trivial action, an
 easy calculation shows that the ring $Q[Z,Z^{-1};\sigma]$  is simple, hence does not admit any
 bounded module and the first  theorem above follows.



For a commutative power series ring in several variables over a
complete discrete valuation ring the Weierstrass preparation
theorem is one key means of proving the unique factorization
property. At least for the completed group algebra (of the
semidirect product $G$) with coefficients in the finite field
$\mathbb{F}_p$ a similar result holds which can hopefully be
extended to the full Iwasawa algebra with \zp-coefficients:

 \begin{bthm}[Theorem \ref{ufd-thm}]
The ring $\fp \kl G \kr=\fp [[X,Y;\sigma,\delta]]$ is an UFD in
the sense of Chatters \cite{cha}.
\end{bthm}

For the precise definition  of a (noncommutative) unique
factorization domain (UFD) we refer the reader to section
\ref{ufd}. To the author's  knowledge this is the first example of
a {\em complete} noncommutative ring of global cohomological
dimension greater than $1,$ which is a UFD.

 \textsc{Acknowledgements.}
 I would  like to thank John Coates most warmly not only for inviting me to  DPMMS,  Cambridge
 University,  but also for his great interest and helpful comments on this project.
 Also I am very  grateful  to Yoshitaka Hachimori for valuable,  stimulating discussions and reading parts of the
 manuscript.  Susan Howson is thanked for helpful suggestions. Finally, I   thank the EU
 Research Training Network ``Arithmetic Algebraic Geometry" for its support and  DPMMS for its hospitality.  \\

{\sc Convention.} In this paper, the ring $R$ is always
associative and with a unit element. When we are talking about
properties related to $R$ like being ``noetherian", an ``ideal", a
``unit", of a certain ``global dimension", etc.\ we always mean
the left {\em and} right property if not otherwise stated. But by
an $R$-module we usually mean left $R$-module  (not a bi-module).
By a {\em local} ring $R$ we mean a ring in which the non-units
form a proper ideal, which is then automatically maximal as left,
right and two-sided ideal. Equivalently, $R$ has both a unique
left and a unique right maximal ideal, which amounts to the same
as the quotient $R/J(R)$ of $R$ by its Jacobson radical being a
skewfield.

\section{Skew power series rings}\label{powerseries}

We begin by recalling the definition of skew power series rings,
which have already been extensively studied in the literature (see
\cite{mc-rob}, cite{li}, \cite{cohn}). Let $R$ be a ring,
$\sigma:R\to R$ a ring endomorphism and $\delta:R\to R$ a
$\sigma$-derivation of $R,$ i.e.\ a group homomorphism satisfying
\[\delta(rs)=\delta(r)s + \sigma ( r) \delta (s)\; \mbox{for all } r,s\in R.\] Then the (formal) skew power
series ring \[R[[X;\sigma,\delta]]\] is defined to be the ring
whose underlying set consists of the usual formal power series
$\sum_{n=0}^\infty r_n X^n,$ with $r_n\in R.$ However, the
multiplication of two such power series is based on the right
multiplication of $X$ by elements of $R$ which is defined by the
formula
\begin{eqnarray}\label{rel} Xr=\sigma (r )X+\delta (r).\end{eqnarray} This clearly implies that, for all
$n\geq 1,$ we have \[X^nr=\sum_{i=0}^n (X^nr)_i X^i\] for certain
elements $(X^nr)_i\in R.$ In other words, we have and in general
the multiplication is given by
\begin{eqnarray*}
  (\sum r_i X^i)(\sum s_jX^j)&=&\sum_n  c_n X^n
\end{eqnarray*}
where  \begin{eqnarray}\label{coeff}
c_n=\sum_{j=0}^n\sum_{i=n-j}^\infty r_i (X^i
s_j)_{n-j}.\end{eqnarray}  To grant that the above series
converges to an element $c_n\in R$ we have to impose further
conditions: Either $\delta =0$ (then $X^nr=\sigma^n(r)X^n,$ i.e.\
the sum $c_n=\sum_{i+j=n} r_j \sigma^i(s_j) $ is finite) or we
assume that $R$ is a complete ring with respect to the $I$-adic
topology, $ I$ some $\sigma$-invariant  ideal, and that it holds
\[\delta(R)\subseteq I,\;\delta(I)\subseteq I^2.\] This implies by induction that
$\delta(I^k)\subseteq (I^{k+1})$ for all $k\geq 0$ ($I^0=R$ by
convention) and the following lemma shows that the multiplication
law in $R[[X;\sigma,\delta]]$ is well defined.

\begin{lem}\label{lemma1}
\begin{enumerate}
\item If $r\in I^k$ then $(X^nr)_i\in I^{k+n-i}$ for all $i,n\geq 0.$
\item If $\sigma$  is an isomorphism satisfying $\sigma(I)=I,$ then for any $n\in\mathbb{N}:$
$r\in I^k\setminus I^{k+1}$ if and only if $(X^nr)_n\in
I^k\setminus I^{k+1}.$ In particular, $r\neq 0$ if and only if
$(X^nr)_n\neq 0.$
\end{enumerate}
\end{lem}

\begin{proof}
All statements follow by induction, the case $n=1$ following from
the defining relation $Xr=\sigma (r )X+\delta r,$ while the
induction step uses
\[(X^{n+1}r)_j=\sigma((X^nr)_{j-1})+\delta((X^nr)_j)\] and $(X^{n}r)_{n}=\sigma^{n}(r).$
\end{proof}

Our interest in this kind of rings stems from applications in
algebraic number theory where the following example naturally
occurs, see section \ref{tate}.

\begin{ex}\label{semi}
Let $G$ be the semi-direct product $G=\Gamma_1 \rtimes \Gamma_2$
of the pro-$p$-groups $\Gamma_1$ and $\Gamma_2,$ both isomorphic
to the additive group of $p$-adic integers $\zp,$ where the action
of $\Gamma_2$ on $\Gamma_1$  is given by a continuous group
homomorphism $\rho:\Gamma_2\to
\mbox{Aut}(\Gamma_1)\cong\mathbb{Z}_p^\ast.$ Then the completed
group algebras $\zp\kl G\kr$ and $\fp\kl G\kr$ of  $G$ with
coefficients in $\zp$ and $\fp$ respectively are isomorphic to the
skew power series rings

\begin{eqnarray*}
\mathbb{Z}_p\kl G\kr
&\cong&\mathbb{Z}_p[[X]][[Y;\sigma,\delta]]=:\mathbb{Z}_p[[X,Y;\sigma,\delta]], \\
\mathbb{F}_p\kl G\kr
&\cong&\mathbb{F}_p[[X]][[Y;\sigma,\delta]]=:\mathbb{F}_p[[X,Y;\sigma,\delta]],
\end{eqnarray*}
where \begin{enumerate}
\item $X:=\gamma_1 -1,\; Y:=\gamma_2 -1$ for some generators $\gamma_i$ of $\Gamma_i,$
$i=1,2.$
\item Setting $R:= \zp[[X]]$ $(\fp[[X]])$ and $\epsilon:=\rho(\gamma_2)$ the ring
automorphism $\sigma$ is  induced by $X\mapsto (X+1)^\epsilon -1,$
while $\delta=\sigma- \mbox{id.}$
\end{enumerate}
Indeed, conferring Lazard 's work \cite{la} (see also \cite{dsms})
the ring $\zp\kl G\kr$ ($\fp\kl G\kr$) has $X^iY^j,$ $0\leq
i,j<\infty,$ as a (complete) $\zp$($\fp$)-basis and it is straight
forward to verify that the relation between $X$ and $Y$ is
described via $\sigma$ and $\delta.$ Furthermore, $R$ is obviously
complete with respect to the topology induced by its maximal ideal
$\m$ which is generated by $X$ and $p.$ Since $\sigma$ is induced
by choosing another generator of $\Gamma_1$ under the isomorphism
$\zp[[X]]\cong \zp\kl \Gamma_1\kr$  we see that $\sigma(\m)=\m.$
Thus we only have to verify that $\delta(\m^k)\subseteq\m^{k+1}$
for $k=0,1:$ Since $\delta$ and $\sigma$ are $\zp$ ($\fp$)- linear
the claim follows immediately from the special case $r=X$
observing also that $\delta(\zp)=0$ (and analogously for $\fp$):

\begin{eqnarray*}
\delta X &=&\sigma(X)-X\\
         &=&(X+1)^\epsilon -1 -X\\
         &=&\sum_{i\geq1} \left(\begin{array}{c}
           \epsilon  \\
           i\
         \end{array}\right) X^i -X\\
         &=&(\epsilon-1)X^k + \mbox{terms of higher degrees}.
\end{eqnarray*}

Since $p|(\epsilon-1)$ we conclude $\delta(X)\subseteq\m^2.$
\end{ex}

\begin{ex}\label{semi-gen}
The previous example generalizes immediately to the following
situation: Let $H$ be a uniform pro-$p$ group on which the group
$\Gamma\cong\zp$ acts via a continuous group homomorphism
$\rho:\Gamma\to \mbox{Aut}(H)$ and assume that the image of $\rho$
is contained in $\Gamma(H^p):=\{\gamma\in Aut(H)|
[\gamma,H]\subseteq H^p\}$ where $[\gamma,h]=h^\gamma h^{-1}$ for
$h\in H.$ Then it is easy to see that  the semidirect product
$G=H\rtimes\Gamma$ is again an uniform pro-$p$ group.  Its Iwasawa
algebras $\Lambda(G):=\zp\kl G\kr$ and $\Omega(G):=\fp\kl G\kr$ of
$G$ with coefficients in $\zp$ and $\fp$  are isomorphic to the
skew power series rings over the Iwasawa algebras $\Lambda(H) $
and $\Omega(H),$ respectively,
\begin{eqnarray*} \Lambda(G)
&\cong&\Lambda(H)[[Y;\sigma,\delta]], \\ \Omega(G)
&\cong&\Omega(H)[[Y;\sigma,\delta]],
\end{eqnarray*}
where $ Y:=\gamma -1$ for some generators $\gamma$ of $\Gamma,$
the ring  automorphism
 $\sigma$ is induced by $h\mapsto h^\gamma$ and again $\delta=\sigma- \mbox{id.}$
 Indeed,
since for any $h\in H$
\begin{eqnarray*}
\delta (h-1) &=&\sigma(h-1)-(h-1)\\
         &=& h^\gamma -h\\
         &=&(g^p-1)h \mbox{ (for some $g\in H   $ by assumption)}\\
         &\in &{\m}^2 \mbox{(see the claim in the proof of \cite[lemma 3.24]{ven1}),}
\end{eqnarray*}
$\delta(\m)\in \m^2$ and the conditions on $ \sigma$ and $\delta$
are fulfilled as above.

We call a uniform pro-$p$-group $G$ {\em polycyclic}     if it
contains a finite chain of subgroups
\[1=G_0\subseteq G_1\subseteq\cdots G_i\subseteq G_{i+1}\subseteq \cdots G_n=G\]  such that
$G_i$ is normal in $G_{i+1}$ with quotient $G_{i+1}/G_i\cong \zp$
for $0\leq i<n.$ Then, as S. Howson suggests, one can identify
$\La(G)$ with the skew power series ring in $n$-variables
\[\Lambda(G)\cong\zp[[X_0,\ldots,X_{n-1};\sigma_1,\delta_1,\ldots,\sigma_n,\delta_n]] \] which is
defined recursively by
\begin{eqnarray*} \Lambda (G_1)&=&\mathbb{Z}_p[[X_0]],\\
\Lambda (G_{i+1})&=&\Lambda (G_i)[[X_i;\sigma_i,\delta_i]],
\end{eqnarray*}
where $\sigma_i$ and $\delta_i$ are defined as above.
\end{ex}

\begin{ex}\label{gl2}
A rather trivial special case of the previous example   is
$R[[X]]$ where $R$ is any (noncommutative) ring,
$\sigma=\mathrm{id}$ and $\delta=0.$ In the Iwasawa theory of
elliptic curves there occur naturally  open subgroups $G\subseteq
GL_2(\zp)$ such that $G=H\times \Gamma$ where $H=SL_2(\zp)\cap G$
and $\Gamma\cong\zp$ is the center of $G.$ Then the Iwasawa
algebra $\La(G)$ can be identified with the power series ring
$\La(H) [[T]]$ in one (commuting) variable with coefficients in
the sub Iwasawa algebra  $\La(H).$
\end{ex}

In case $\delta=0$ the following properties of $A=R[[X;\sigma]]$
reflect the corresponding ones of $R$ (\cite[thm. 1.4.5, thm.
7.5.3]{mc-rob}, \cite[Chap. II Cor. 3.2.12, Chap. III Thm.
3.4.6]{li}):

\begin{prop}
Let $A=R[[X;\sigma]]$. Then the following holds:
\begin{enumerate}
\item If $\sigma$ is injective and $R$ is an integral domain, then $A$ is an integral domain.
\item If $\sigma$ is injective and $R$ a division ring, then $A$ is a principal left ideal
domain (Note that in this case $\delta=0$ because of the above
conditions).
\item If $\sigma$ is an automorphism and $R$ is left (or right) Noetherian, then $A$ is left
(respectively right) Noetherian. In this case, the left (or right)
global dimensions of $A$ and $R$ are related by
\[\mathrm{gl}A=\mathrm{gl}R+1.\]
\item If $\sigma$ is an automorphism and $R$ is a maximal Noetherian order, then $A$ is a maximal
 Noetherian order.
 \item If $\sigma$ is an automorphism and $R$ an Auslander regular ring, then $A$ is an Auslander regular ring.
\end{enumerate}
\end{prop}

In order to prove similar results in the general case we first
have to introduce some filtrations on $A= R[[X;\sigma,\delta]].$
We assume for the rest of this section that $R$ is a local ring
with maximal ideal $\m$ which is complete and separated with
respect to its $\m$-adic topology. Then we set $F_0A:=A$ and
define, for $n\geq 1,$   $F_nA$ to be the set of all
$f=\sum_{k=0}^\infty c_kX^k$ such that $c_k\in\m^n$ for all $k\geq
0.$
Finally, we set
\[G_0:= A,\;\; G_n:=\sum_{i=0}^n (F_{n-i}A)X^i,\;\ n\geq 1,\] where  the sum is a sum of left
By definition, both filtrations $G_i$ and $F_i$ are exhaustive and
$F_iA\subseteq F_jA$ respectively $G_iA\subseteq G_jA$ holds for
all $i\geq j.$ The following lemma shows that they are indeed
filtrations of rings:

\begin{lem}
For all $n,$  $m\geq 0,$ it holds: $F_nA \cdot F_m A\subseteq
F_{n+m} A$ and $G_nA \cdot G_m A\subseteq G_{n+m} A.$
\end{lem}

\begin{proof} The assertion for the filtration $F_i$ follows immediately from the formula
\ref{coeff}, lemma \ref{lemma1} and by the fact  that $\m^k$ is
closed with respect to the $\m$-adic topology for all $k$ (see
\cite[Ch. I 3.1(c)]{li}). Concerning the second filtration it
suffices by distributivity to prove that
\[(F_{n-i}A)X^i(F_{m-j}A)X^j\subseteq \sum_{k=0}^{n+m}(F_{n+m-k}A)X^k\] for all $0\leq i\leq
n$ and $0\leq j\leq m$ holds. Using the multiplicativity of $F_i$
just proved this will follow immediately once we have shown the
relation \[X^i(F_kA)\subseteq \sum_{j=0}^i (F_{k+i-j}A)X^j.\] So
let $f=\sum f_i X^i$ in $F_kA.$ The defining relation \ref{rel}
implies that $Xf$ decomposes into the sum
\[Xf=\sum_{i=0}^\infty\delta(f_i)X^i + (\sum_{i=1}^\infty
\sigma(f_{i-1})X^{i-1})X,\] the first summand of which lies in
$F_{k+1}A$ while the second one belongs to $(F_kA)X$ by lemma
\ref{lemma1}. The proof is accomplished by induction.
\end{proof}

Note that the two-sided ideals $G_kA$ (apply the lemma with $n$ or
$m$ equal to $0$) can be described as follows: \begin{eqnarray}
\label{prod-des} G_kA=(\prod_{i=k}^\infty RX^i)\times\m
X^{k-1}\times \cdots \times\m^k X^0=\prod_{i=0}^\infty \m^{k-i}
X^i\end{eqnarray} for all $k\geq 0,$ where $\m^l:=R$ for negative
integers $l$ and where $A$ is identified with
$A=\prod_{i=0}^\infty R X^i.$

\begin{lem}
The topology on $A$ induced by the filtration $G_iA$ coincides
with the product topology of $A\cong\prod_{i=0}^\infty R,$ where
$R$ is endowed with its $\m$-adic topology.
\end{lem}
\begin{proof}
With respect to the product topology a neighbourhood basis of
$0\in A$ consists of
\[A_k:=\prod_{i=0}^{k-1}\m^kX^i\times\prod_{i=k}^\infty R X^i,\; k\geq 0.\] It is easy to see that the
$B_k$ and $G_i$ are cofinal to each other and hence induce the
same topology on $A.$
\end{proof}

Since the product of complete Hausdorff topological spaces is
again a complete Hausdorff topological space we obtain

\begin{cor}
The ring $A$ is a  complete Hausdorff topological ring with
respect to its product topology respectively the topology which is
induced by the filtration $G_i.$
\end{cor}

\begin{prop}
There is a canonical isomorphism of graded rings
\[gr_GA\cong (gr_\m R)[\bar{X};\bar{\sigma}],\] where the latter ring bears the mixed grading
with respect to the grading of $gr_\m R$ and the "$\bar{X}$-adic"
grading. Here $\bar{X}$ means the image of $X$ and $\bar{\sigma}$
is induced by $\sigma.$
\end{prop}

\begin{proof}
The above description \ref{prod-des} implies immediately that
\[G_k/G_{k+1}\cong\prod_{i=0}^k\m^i/\m^{i+1} \bar{X}^{k-i}\] holds. Thus
\begin{eqnarray*}
\bigoplus_{k=0}^\infty G_k/G_{k+1}&\cong&\bigoplus_{i,l\geq 0} \m^i/\m^{i+1} \bar{X}^l\\
      &\cong&\bigoplus_{i\geq 0} gr_\m R \bar{X}^i\\
      &\cong& (gr_\m R)[\bar{X};\bar{\sigma}].
\end{eqnarray*}
\end{proof}

As a corollary it turns out that $A=R[[X;\sigma,\delta]]$ reflects
the properties of $gr_\m R:$

\begin{cor}\label{grR}
If $gr_\m R$ is

\begin{enumerate} \item integral and $\bar{\sigma}$ is injective, then $A$ is
integral. \item left (or right) Noetherian and $\bar{\sigma}$ is
an automorphism, then $A$ is left (respectively right) Noetherian.
If, in addition, the left (or right) global dimension of $gr_\m R$
is finite, then so is the left (respectively right) global
dimension of $A$ and it holds
\[\mathrm{gl}A\leq\mathrm{gl}gr_\m R+1.\]
\item a Noetherian maximal order, $R$ complete and  $\bar{\sigma}$ is an automorphism,  then $A$ is a maximal Noetherian order.
\item  Noetherian, $R$ complete and $\bar{\sigma}$ is an automorphism, then $A$ is
 a Zariski-ring.
\item Auslander regular, $R$ complete and $\bar{\sigma}$ is an automorphism,
 then $A$ is Auslander regular.
\end{enumerate}
\end{cor}

\begin{proof}
 Assuming conditions (i) or (ii), \cite[thm. 1.2.9]{mc-rob} tells that \linebreak
 $ gr_GA\cong (gr_\m R)[\bar{X};\bar{\sigma}]$ is integral, respectively Noetherian. This
 implies the same properties for the ring $A$ by \cite[thm 1.6.7, thm. 1.6.9]{mc-rob}. The
 statement concerning the global dimensions follow similarly by \cite[thm. 7.5.3 (iii), cor.
 7.6.18]{mc-rob}.
 Now let us assume (iii). Then  $gr_GA\cong (gr_\m R)[\bar{X};\bar{\sigma}]$ is a maximal
 order according to \cite[Chap. II Cor. 3.2.12]{li}, thus $A$ is a maximal order by
 \cite[lem 2.6]{co-sch-su}. Under the condition (iv) $A$ is complete with
 respect to the filtration $G_i.$ Since $gr_G A$ is Noetherian, $A$ is a Zariski ring by
 \cite[Chap. II thm. 2.1.2]{li}. To prove (v) it suffices to show that $gr_GA$ is  Auslander
 regular (\cite[Chap. III thm. 2.2.5]{li}), which is contained in \cite[Chap. III thm. 3.4.6
 (1)]{li}.
\end{proof}

We close this section proving that under our current conditions on
$R$ the skew power series ring   $A$ is a  local ring.

\begin{prop}\label{local}
An element $f=\sum f_i X^i\in A$ is a unit (in $A$) if and only if
the constant term $f_0$ is a unit in $R.$ In particular, $A$ is a
local ring.
\end{prop}

\begin{proof}
The ``only if" part follows by considering the canonical
surjective map $A\to A/G_1A\cong R/\m.$ Now assume  without loss
of generality that $f_0=1$ (otherwise replace $f$ by $f_0^{-1}f$ )
and set $h:=f_0-f\in G_1A.$ Then, due to completeness of $A$ with
respect to $G_iA,$ the ``geometric series" $\sum_{i\geq 0} h^i$
converges to an element $g\in A,$ which is a left and right
inverse of $f=1-h.$ That means that the non-units consist
precisely of the elements of $\M=G_1A.$
\end{proof}

\section{The Weierstrass preparation theorem}

Let $R$ be a (not necessarily commutative) local ring with
maximal (left) ideal $\m$  and suppose that $R$ is separated and
complete with respect to its $\m$-adic topology. As usual we
denote the residue class skewfield $R/\m$ by $k.$ As before we
assume that $\sigma:R\to R$ is a ring-automorphism satisfying
$\sigma(\m)=\m,$ i.e.\ $\sigma$ induces also a ring automorphism
$\sigma:k\to k,$ while $\delta:R\to R$ is a $\sigma$-derivation
such that $\delta(R)\subseteq\m,\; \delta(\m)\subseteq\m^2.$ In
particular, $\delta$ operates trivially on $k.$ Therefore there is
a canonical surjective reduction map:

\[\bar{ } :R[[X;\sigma,\delta]]\to k[[X;\sigma]].\]

In example \ref{semi} $\sigma$ will turn  out to be the identity
on $k.$ For any $f=\sum a_i X^i\in A:=R[[X;\sigma,\delta]]$ the
reduced order $\mbox{ord}^{red}(f)$  of $f$ is defined to be the
order of the reduced power series $\overline{f},$ i.e.\
 \[\mbox{ord}^{red}(f)=\min\{i|a_i\in R^\ast\}.\]

 \begin{thm}\label{WP}
 Let $f\in A$ be a power series with finite reduced order $s:=\mbox{ord}^{red}(f)<\infty.$
  Then $A$ is the direct sum
 \[A=Af\oplus\bigoplus_{i=0}^{s-1}RX^i\]
 of the $R$-modules $N:=Af$ and $M:=\bigoplus_{i=0}^{s-1}RX^i.$
 \end{thm}

 \begin{proof}
 (compare to Bourbaki's proof in the commutative case \cite[{\S}3 no.8]{bourbaki}) We first prove
 that $N\cap M=0:$

 For any element of this intersection one has an equality
 \begin{eqnarray*}
r_0+r_1X+\ldots+r_{s-1}X^{s-1}&=& (\sum b_j X^j)f\\
                              &=&\sum c_nX^n.
 \end{eqnarray*}

 Comparing the coefficients we conclude that for $k\geq 0$
 \begin{eqnarray*}
 0\stackrel{!}{=}c_{s+k}&=&\sum_{\begin{array}{c}
   0\leq j\leq s+k \\
   i \
 \end{array}}b_i (X^ia_j)_{s+k-j}\\
 &=& b_k(X^ka_s)_k + \sum_{\begin{array}{c}
   0\leq j\leq s+k \\
   i \\
   (i,j)\neq (k,s)
 \end{array}}b_i (X^ia_j)_{s+k-j}.
 \end{eqnarray*}
 To prove that $b_j=0, j\geq 0,$ it suffices to see that $b_k\in\m^n$ for all $k$ and $n.$ By double
 induction, let us assume that $b_i\in\m^{n-1}$ for all $i$ and $b_i\in\m^n$ for $i<k.$ Since
 $a_s\in R^\ast,$ so is $(X^ka_s)_k$ by lemma \ref{lemma1} and the following claim implies
 that $b_k\in\m^n.$\\

 {\em Claim:} $b_i(X^ia_j)_{s+k-j}\in\m^n$ for all $0\leq j\leq s+k,\;(i,j)\neq (k,s).$

 For $i<k$ the claim is obvious because of $b_i\in\m^n$ then. So let us assume $i\geq k,
 (i,j)\neq (k,s).$ We have to deal with three cases separately:

 \begin{enumerate}
 \item If $i+j\leq s+k$ then $(X^ia_j)_{s+k-j}=0.$
 \item If $i+j=s+k$ then $j=s+k-i<s$ (note that $i=k$ would imply that $j=s$), i.e.\
 $a_j\in\m$ by the definition of $s$
  and the claim follows from lemma
 \ref{lemma1}.
 \item If $i+j>s+k$ then $i-(s+k-j)\geq 1,$ i.e.\ the claim follows from the same lemma.
 \end{enumerate}

 We still have to prove that $N+M=A.$ Since $a_s\in R^\ast$ and $A$ is local we can write
 \[f=\sum_{i=0}^{s-1}a_iX^i +gX^s\]
 with $g\in A^\ast.$ Setting $G:=g^{-1}$ and
 \begin{eqnarray*}
 h:=-G(\sum_{i=0}^{s-1}a_iX^i)&=&-G(f-gX^s)=X^s-Gf
 \end{eqnarray*}
 we claim that all coefficients of $h$ are in $\m.$ Indeed,
 \[h_n=\sum_{j\leq\min\{n,s-1\}} G_i (X^ia_j)_{n-j}\in \m\]
  using lemma \ref{lemma1} because $a_j\in \m$ for $j<s.$

  Now let $r$ be an arbitrary element in $A.$ Recursively we can find a series $q^{(n)}=\sum_{l=0}^\infty q^{(n)}_l X^l,
  n\in\mathbb{N},$ of elements in $A$ such that
  \[q^{(0)}X^s\equiv r \mbox{ mod }M\]
  and
  \[q^{(n+1)}X^s\equiv q^{(n)}h \mbox{ mod } M \mbox{ for }n\geq0\]
 (Just cut off the terms of degrees less than $s$ with respect to $X$ and then factor out
 $X^s$ on the right).
 Since $q^{(n+1)}_l=(q^{(n)}h)_{s+l}$ it is easily seen that $q^{(n+1)}_l\in\m^{n+1}$ for all $l$ using the
 fact that the coefficients of $h$ lie all in $\m.$ Hence, due to the completeness of $R,$
  $\sum q^{(n)}$ converges to an element $q\in A.$ The relations
  \[(q^{(0)}+\dots +q^{(n)})X^s\equiv r+ (q^{(0)}+\dots +q^{(n-1)})h \mbox{ mod } M\] imply that
\begin{eqnarray*}
 r&\equiv& q(X^s-h)\mbox{ mod } M\\
 &=&(qG)f\in Af.
\end{eqnarray*}
 \end{proof}

In analogy with the commutative situation we call a monic
polynomial $F=X^s+a_{s-1}+\cdots +a_1 X + a_0 \in
R[X;\sigma,\delta]\subseteq R[[X;\sigma,\delta]]$ {\em
distinguished} or {\em Weierstrass polynomial} (of degree $s$) if
$a_i\in\m$ for all $0\leq i\leq s-i.$ Here, $R[X;\sigma,\delta]$
denotes the skew polynomial ring in one variable over $R$ with
respect to $\sigma$ and $\delta.$

\begin{cor}\label{cor-WP}
Under the assumptions of the theorem $f$ can be expressed uniquely
as the product of an unit $\epsilon$ of $A$ and a distinguished
polynomial $F\in A:$ \[f=\epsilon F.\]
\end{cor}

\begin{proof}
According to the theorem the element $X^s\in A$ can be uniquely
written as
 \[X^s=vf-H\]
where $H=\sum_{i=0}^{s-1} g_i X^i$ and $v\in A.$ Reducing the
coefficients $\mbox{mod }\m$ we obtain \[X^s=\bar{v} \bar{f} -
\sum_{i=0}^{s-1}\bar{ g_i} X^i .\] Since
$\mbox{ord}(\bar{f})=\mbox{ord}^{red}(f)=s,$ i.e.\
$\mbox{ord}(\bar{v}\bar{f})=\mbox{ord}(\bar{v})+\mbox{ord}(\bar{f})\geq
s,$ the coefficients $\bar{g_i}$ must all vanish. Therefore
\[F:=X^s+H\] is a distinguished polynomial and
\[X^s=\bar{v}\bar{f},\]
i.e.\ $\mbox{ord}(\bar{v})= \mbox{ord}(X^s)-\mbox{ord}(\bar{f})
=0.$ Hence, $v$ is a unit in $A$ and $f=\epsilon F$ with
$\epsilon:=v^{-1}.$
\end{proof}

Independently, D. Burns and C. Greither \cite{burns-greither} have
proved another generalized version of the Weierstrass preparation
theorem for power series rings in one commuting variable over
noncommutative (not necessarily local) rings.

\begin{rem}
Let $J$ be a left ideal of $A=R[[X;\sigma,\delta]].$ Then $A/J$ is
finitely generated as $R$-module if and only if $J$ contains a
Weierstrass polynomial. Indeed, if all elements of $J$ reduce to
zero in $\kappa:=k[[X;\sigma]],$ then $M/\M M\cong\kappa$ where
$\M:=\ker(\bar{ }:A\to \kappa).$ But since there is some
surjection $R^n\to M$ and $\m\subseteq\M$ the module $M/\M M$ is a
finitely generated $k$-module, a contradiction. The other
implication is a direct consequence of theorem \ref{WP}.

If, moreover, $J=Af$ is principal, then $A/J$ is a finitely
generated $R$-module if and only if $\mbox{ord}^{red}(F)$ is
finite, i.e.\ if and only if $Af=AF$ can also be generated by some
Weierstrass polynomial $F.$
\end{rem}

Of course, there exist right versions for all statements of this
section. Indeed, since $\sigma$ is assumed to be an automorphism,
one gets a natural ring isomorphism

\[R[[X;\sigma,\delta]\cong[[X;\sigma',\delta']]R,\]
where  $\sigma'=\sigma^{-1},$  $\delta'=-\delta\circ\sigma^{-1}$
and in the latter ring the coefficients are written on the right
side of the variable $X.$ Moreover, the reduced order is invariant
under this isomorphism, i.e.\ the notion of an distinguished
polynomial and its degree is independent of the representation as
left or right power series.

\section{Faithful modules}\label{faithfulmod}

In this section we calculate the global annihilator ideals of
certain modules over the Iwasawa algebra of some $p$-adic Lie
group. First we concentrate on the example \ref{semi}, i.e.\
$A=R[[Y;\sigma,\delta]]$ with $R=\zp[[X]]$ or $R=\fp[[X]].$ Recall
that in this case $\delta$ is just given as $\sigma -\mbox{id}.$
To exclude the case where $G$ respectively $A$ are commutative, we
assume that
 $\epsilon\neq 1.$

\begin{rem}\label{rem-ann}
\begin{enumerate}
\item There is a canonical isomorphism
\[R[Y;\sigma,\delta]\cong R[Z_;\sigma],\;r\mapsto r,\;Y\mapsto Z-1\] of skew polynomial rings, because
\begin{eqnarray*}
Zr&=&(Y+1)r\\
  &=&Yr +r\\
  &=&\sigma(r)Y +\delta r+r\\
  &=&\sigma(r) Y + \sigma(r)-r+r\\
  &=&\sigma(r)(Y+1)=\sigma(r) Z.
\end{eqnarray*}
\item Any  ideal $I$ of $A$ which contains a polynomial $0\neq b=\sum_{i=0}^sb_iZ^i$
also contains a non-zero element $r\in R.$ First note that
$Z^i\gamma=\sigma^i(\gamma)Z^i$ for $\gamma:=X+1\in R^\ast,$ i.e.\
$
\sigma^j(\gamma)Z^i-Z^i\gamma=(\sigma^j(\gamma)-\sigma^i(\gamma))Z^i$
the latter being zero if and only if $i=j$ because $\sigma$
operates without fixpoint on $\Gamma_1.$   Since $I$ is two-sided
it also contains the element
\begin{eqnarray*}
\sigma^s(\gamma)b-b\gamma&=&\sum_{i=0}^s b_i(\sigma^s(\gamma)Z^i-Z^i\gamma)\\
                         &=&\sum_{i=0}^{s-1}b_i (\sigma^s(\gamma)-\sigma^i(\gamma))Z^i
\end{eqnarray*}
which is nonzero whenever one of the $b_i, 0\leq i\leq s-1,$ is
nonzero because $R$ is integral. Proceeding recursively one
concludes that $I$ contains an element of the form $rZ^i$ for some
$0\neq r\in R$ and $i\in \mathbb{N}.$ But then it contains also
$r$ itself because $Z=Y+1$ is a unit in $A.$
\item The same argument proves that the skew Laurent polynomial ring \linebreak $Q[Z,Z^{-1};\sigma]$ over the
maximal ring of quotients $Q$ of $R$ is simple (note that the ring
automorphism $\sigma$ extends uniquely to $Q.$)
\end{enumerate}
\end{rem}

\begin{prop}\label{faithful}
\begin{enumerate}
\item Let $R=\fp[[X]].$ If $M$ is an $A$-module which is not  torsion as
$R$-module, then its annihilator ideal over $A$ vanishes:
\[Ann_A(M)=0.\]
\item Now let $R=\zp[[X]].$ If $M$ is an $A$-module without  $p$-torsion and
such that the module $M/p$ is not  torsion as
$R/p=\fp[[X]]$-module (for example if $M$ is a free $R$-module),
then $Ann_A(M)=0.$
\end{enumerate}
\end{prop}

Recall that an $A$-module $M$ is called {\em faithful} if the
annihilator ideal $Ann_A(M)$ is zero and otherwise $M$ is called
{\em bounded}. Before proving the proposition we want to draw an
immediate conclusion from it in the case $R=\zp[[X]].$ If $M$ is
an $A$-module which is finitely generated as $R$-module and of
strictly positive $R$-rank $\rk_R M>0,$ then the conditions of
(ii) in the proposition are satisfied for $N:=M/\tor_\zp M.$
Indeed,  by \cite[cor. 1.11]{howson2000} it holds
$\rk_{R/p}N/pN=\rk_{R/p}({_pN)}+\rk_R N =\rk_RN\neq 0$ in this
case, where $_pN$ denotes the kernel of multiplication by $p.$
Thus

\begin{cor}\label{cor-faithful}
Every   $A$-module $M$ which is finitely generated as $R$-module
and of strictly positive $R$-rank $\rk_R M>0$ is faithful.
\end{cor}


\begin{proof} [Proof (of prop. \ref{faithful})]
By the following lemma the second statement is a consequence of
the first one, which can be proven as follows. First note that the
$R$-torsion submodule $M_{R-\text{tor}}$ of $M$ is an $A$-module,
because $X$ is a normal element of $A$ (see below). Thus $M$ has a
quotient $N$ without $R$-torsion and since $Ann_A(M)\subseteq
Ann_A(N)$ we may assume without loss of generality that $M$ itself
is a torsionfree $R$-module. Now,  if $Ann_A(M)$ contains a
non-trivial element $b\in A$ this can be written as
$X^n\widetilde{b}$ with $\mbox{ord}^{red}(\widetilde{b})<\infty$
because $R$ is a principal ideal domain all ideals of which are
generated by some power of $X.$ Due to the Weierstrass preparation
theorem $\widetilde{b}$ is a Weierstrass polynomial in the
variable $Y$ up to a unit, i.e.\ we may assume that $b$ is already
a Weierstrass polynomial because neither the unit nor $X^n$ can
annihilate $M$ by assumption. By the above remark then the
annihilator ideal contains a nonzero element $r\in R$ which
contradicts the assumption that $M$ is a torsionfree $R$-module.
\end{proof}

\begin{rem}
By similar arguments one can show that every nonzero two-sided
ideal $J$ of $A=\fp[[X,Y;\sigma,\delta]]$ contains some power of
$X.$ It follows immediately that for every two-sided ideal $J$ of
$A=\zp[[X,Y;\sigma,\delta]]$ the cyclic module $A/J$ is finitely
generated over $\zp[[Y]]\cong\La(\Gamma)$ (for any lift
$\Gamma\subseteq G$). In particular, by an easy induction argument
on the number of generators, each finitely generated bounded
$A$-module is finitely generated over $\La(\Gamma).$
\end{rem}

\begin{lem}
Let $\Lambda:=\zp\kl G \kr$ and $\Omega:=\fp \kl G\kr\cong
\Lambda/p$ the completed group algebras of some pro-finite group
$G$ with coefficients in $\zp$ and $\fp$ respectively and assume
that both rings are Noetherian integral domains. If  $M$ is an
$A$-module without $p$-torsion then the following holds

\[ Ann_\Lambda(M)\neq 0\; \Rightarrow Ann_\Omega(M/p)\neq 0.\]
\end{lem}

\begin{proof}
If $f$ is any nonzero element in $Ann_\Lambda(M)$ and $f=p^nf'$ is
the unique factorization such that $f'$ is not divisible by $p$
then $f'$ is also an element in the annihilator of $M$ because the
latter module does not have any $p$-torsion by assumption. Thus
the image of $f'$ in $\Omega$ is obviously a non-zero element of
$Ann_\Omega(M/p).$
\end{proof}

In the following we study the behaviour of the annihilator ideal
under induction. So let $H$ be a closed subgroup of a $p$-adic Lie
group $G.$ If $M$ denotes a finitely generated
$\Lambda(H)$-module, then we write
$\Ind^H_GM:=\Lambda(G)\otimes_{\Lambda(H)} M$ for the
$\Lambda(G)$-module which is induced from $M.$

Let us assume for a moment that $H=(h_1,\ldots,h_r)$ is a uniform
pro-$p$-group  with minimal system of generators $h_1,\ldots,h_r.$
Then, using the isomorphism $\Lambda(H)\cong\prod_{(l_1,\ldots,
l_r)\in\mathbb{N}^r} \zp (h_1-1)^{l_1}\cdots (h_r-1)^{l_r},$ it is
easily seen that the $\zp$-linear span $\sum\zp h_n$ is dense in
$\Lambda(H)$ where $h_n:=h_1^{\nu^n_1}\cdots h_r^{\nu^n_r}$ for
some bijection $\mathbb{N}\cong\mathbb{N}^r,$ $n\to
\nu^n=(\nu^n_1,\ldots, \nu^n_r).$  Since every $p$-adic analytic
group possesses an uniform subgroup of finite index it follows
immediately that there exists (for arbitrary $H$ as above) a
countable set $\mathcal{H}=\{h_n\in H |n\in\mathbb{N}\}$ such that
its $\zp$-linear span is dense in $\Lambda(H).$ For a left ideal
$I$ and an element $r$ of a ring  $\Lambda$ we will denote by
$(I:r):=\{ \lambda\in\Lambda| \lambda r\in I\}$ the annihilator of
$r+I\in \Lambda/I.$

\begin{lem}\label{ann-intersection}
Let $I$ be a left ideal of $\Lambda:=\Lambda(H)$ and set
$M:=\Lambda/I.$ Then, for any choice of  $\mathcal{H},$ the
annihilator ideal equals the countable intersection
\[Ann_\Lambda(M)=\bigcap_{n=0}^\infty I_n\] of  the left ideals $I_n:=Ih_n^{-1}=(I:h_n).$
\end{lem}

\begin{proof}
In general, it holds
\begin{eqnarray*}
Ann_\Lambda(M)&=&\bigcap_{m\in M} Ann_\Lambda (m) \\
              &=&\bigcap_{r \in \Lambda} (I:r).
\end{eqnarray*}
Thus $Ann_\Lambda(M)$ is contained in $ \bigcap_{n=0}^\infty I_n.$
On the other hand any element $\lambda$ of the latter intersection
annihilates $r+I\in M$ for any $r$ in the $\zp$-linear span of
$\mathcal{H}.$ Since this span is dense in $\Lambda$ the statement
follows by the continuity of the action of $\Lambda$ on $M.$
\end{proof}

In order to compare $Ann_{\Lambda(G)}(\Ind^H_GM)$ with
${\Lambda(G)}Ann_{\Lambda(H)}(M)$ the following lemma will be
crucial. Note that $\Lambda(G)$ is a flat $\La(H)$-module
(\cite[proof of lem 5.5]{ochi-ven}).

\begin{lem}\label{intersection}
Let $I_i,$ $i\in\mathbb{N},$ be a set of left ideals in
$\Lambda(H).$ Then the following holds:
\[\Lambda(G) \bigcap_i I_i = \bigcap _i  \Lambda(G)I_i.\] A similar statement holds for an
arbitrary filtered index set $T$ if the ideals form a descending
chain $I_s\supseteqq I_t$ for $t\geq s.$
\end{lem}

\begin{proof}
Due to the flatness of $\Lambda(G)$ the functor
$\Lambda(G)\otimes_{\Lambda(H)}-$ interchanges with finite
intersections. Hence we can replace the ideals $I_i$ by the ideals
$J_i:=\bigcap_{n=0}^i I_n$ to come into the situation of the
second statement, in which the intersection can be expressed as a
projective limit: $\bigcap_T I_t=\projlim T I_t.$ Since inverse
limits interchange with the completed tensor product
$\La(G)\widehat{\otimes}_{\Lambda(H)} -$ and the latter coincides
with the usual tensor product if applied to finitely generated
$\Lambda(H)$-modules (see \cite{brumer}) we obtain

\begin{eqnarray*}
\Lambda(G) \bigcap_T I_t &=&\Lambda(G)\otimes_{\Lambda(H)} \bigcap_T I_t\\
                          &=&\Lambda(G)\widehat{\otimes}_{\Lambda(H)} \bigcap_T I_t\\
                          &=&\Lambda(G)\widehat{\otimes}_{\Lambda(H)} \projlim T I_t\\
                          &=&\projlim T  \Lambda(G)\widehat{\otimes}_{\Lambda(H)}  I_t\\
                          &=&\projlim T  \Lambda(G)\otimes_{\Lambda(H)}  I_t\\
                          &=&\projlim T  \Lambda(G) I_t.
\end{eqnarray*}
\end{proof}

\begin{prop}
Let $M$ be a finitely generated $\Lambda(H)$-module. Then the
global annihilator ideal of $\Ind^H_G M$ is contained in ideal
generated by the global annihilator of $M:$

\[  Ann_{\Lambda(G)}(\Ind^H_G M)\subseteq\Lambda(G)Ann_{\Lambda(H)}(M).\] In particular, if $M$ is
a faithful $\Lambda(H)$-module, then $\Ind^H_G M $ is a faithful
$\Lambda(G)$-module.
\end{prop}

\begin{proof}
Since a finitely generated $\Lambda(H)$-module $M$ is the sum of
finitely many cyclic submodules $M_i$ and
$Ann_{{\Lambda(H)}}(M)=\bigcap_i Ann_{\Lambda(H)}(M_i)$  this
statement is easily reduced to the case of an cyclic module
$M=\Lambda(H)/I.$  Then $\Ind^H_G M$ is isomorphic to
$\La(G)/\La(G)I$ and it holds

\[Ann_{\La(G)}( \La(G)/\La(G)I )\subseteq \bigcap_{\mathcal{H}} \Lambda(G)Ih^{-1}=\La(G) Ann_{\La(H)}
(M)\] by the previous lemmata
\end{proof}

This result should be compared to a theorem of Harris
\cite{harris-bounded} which tells that the module $\Ind^H_G \zp$
is bounded whenever $2\dim H>\dim G.$ Using this proposition we
can produce a series of faithful $\Lambda(G)$-modules where $G$ is
either an appropriate pro-$p$-subgroup of $SL_2(\zp)$ or
$GL_2(\zp)$  (and thus $SL_n(\zp)$ or $GL_n(\zp)$). Indeed, the
closed subgroup $H$ of $GL_2(\zp)$ generated by the matrices $
t=\left(\begin{array}{cc}
  a & 0 \\
  0 & a^{-1}
\end{array} \right)$ and $c=\left(\begin{array}{cc}
  1 & b \\
  0 & 1
\end{array} \right),$ where e.g. $a=1+p$ and $b=1,$  i.e.\ more or less ``the" (pro-$p$-) Borel subgroup of $SL_2(\zp),$ is
isomorphic to the semidirect product $H\cong C\rtimes T$ of
$C=c^\zp=\left(\begin{array}{cc}
  1 & \mathbb{Z}_p \\
  0 & 1
\end{array} \right)$ and $T=t^\zp.$  Hence all the induced modules $\Ind^H_G(\Lambda(H)/\La(H)F)$ where $F$
denotes a distinguished polynomial in $\La(H)\cong
\zp[[X,Y;\sigma,\delta]]$ are faithful. We should mention that the
faithfulness of  $\Ind^H_G(\Lambda(H)/\La(H)(t-1))\cong
\Ind^T_G\zp$ was proved by Greenberg (private communication) using
``$p$-adic harmonic analysis."

\section{Modules up to pseudo-isomorphism}
In this section let $R$ denote a noetherian local ring with
maximal ideal $\m_R$ which is {\em compact} with respect to its
$\m_R$-adic topology. Furthermore, let  $\sigma$ be a ring
automorphism and $\delta$ a $\sigma$-derivation as in section
\ref{powerseries}. We  write $A$ for the (compact)  ring
$R[[Y;\sigma,\delta]]$ and $B$ for the (non-compact) ring
$R[Y;\sigma,\delta].$ If $\m_A$ denotes the unique maximal ideal
of $A,$ then $\m_B:=B\cap\m_A$ is a maximal  ideal of $B$ and they
are the kernels of the canonical $R$-algebra homomorphisms
\begin{eqnarray*}
{\m_A}&=&\ker(A=R[[Y;\sigma,\delta]]\to R/\m_R),\\ {\m_B}&=&
\ker(B=R[Y;\sigma,\delta]\to R/\m_R),
\end{eqnarray*}
which are both induced by $Y\mapsto 0.$

\begin{lem}\label{m-adic} With the above notation it holds that
\begin{enumerate}
\item $\m_A=\m_RA +YA=A\m_R+AY$ and $\m_B=\m_RB+YB=B\m_R+BY.$
\item For any $A$-module $M,$ the submodules $\m_A^iM=\m_B^iM$ coincide as $B$-modules. In
particular, 
\begin{eqnarray*}
A&=&\projlim{i}A/\m_A^i=\projlim{i}B/\m_B^i \mbox{ and }\\ M&=&
\projlim{i}M/\m_A^iM=\projlim{i}M/\m_B^iM.
\end{eqnarray*}
 \end{enumerate}
\end{lem}

\begin{proof}
Due to compactness and the fact that $\m_R$ is finitely generated
over $R$, it is easy to see that $\m_RA=F_1A$ and thus
$\m_RA+TA=\m_RA+AT=G_1A=\m_A.$ The proof for $B$  is similar, but
easier, while the ``left-versions" follow by symmetry. The second
item follows immediately from the first one.
\end{proof}

By $A$-mod and $B$-mod we denote the category of finitely
generated $A$- and $B$-modules, respectively. For the full
subcategories consisting of modules which are even finitely
generated over $R$ we write $\mbox{$A$-mod}^R$ and
$\mbox{$B$-mod}^R,$ respectively. Similarly, we write mod-$A,$
$^R\mbox{mod-$A$},$ etc.\ for the corresponding categories of
right modules.

\begin{lem}\label{A-sub}
Let $M$ be in $\mbox{$A$-mod}^R.$ Then any $B$-submodule
$N\subseteq M$  is an $A$-submodule and any  $B$-quotient  of $M$
is an $A$-quotient.
\end{lem}

\begin{proof}
The first statement obviously implies the second one. Since all
rings under consideration are noetherian we may assume that
$N=Bm\subseteq M$ is a cyclic $B$-submodule. Then the
$A$-submodule $Am\subseteq M$ is finitely generated over $R$ and
thus the left ideal $Ann_A(m)$ contains a distinguished
polynomial. By the Weierstrass preparation theorem,
$A/Ann_A(m)\cong B/Ann_B(m)$ and the statement for the submodules
follows.
\end{proof}

Recall that the grade $j_A(M)$ is the minimal integer $i$ such
that
\[\E_A^iM:=\Ext_A^i(M,A)\] does not vanish.

\begin{lem}\label{pseudo-torsionI} Let $M$ be in $\mbox{$A$-mod}^R.$ Then the following holds:
\begin{enumerate}
\item  There is a natural identification of $A$-modules
\[M=A\otimes_B M,\] where $M$ is considered as $B$-module via restriction. In particular,
$A\otimes_B -$ is faithfully flat on the full subcategory of
$\mbox{$B$-mod}^R$ whose objects are  restrictions of
$\mbox{$A$-mod}^R.$
\item  For all  $i$ there is a natural isomorphism of right $A$-modules
\[\E^i_B(M)\otimes_B A\cong\E^i_A(M).\] Moreover, $\E^i_A(M)$ and $\E^i_B(M)$ are in $^R\mbox{mod-$A$}$
and $^R\mbox{mod-$B$},$ respectively, and  $\E^i_B(M)$ bears a
natural right $A$-module structure which extends the right
$B$-module structure.
\item  The grades of $M$ as $A$- and $B$-module coincide: \[j_A(M)=j_B(M).\] In particular, $M$ is a
pseudo-null $A$-module if and only if it is a pseudo-null
$B$-module.
\end{enumerate}
\end{lem}

\begin{proof}
By \cite[Ch. II, prop. 1.8 (3)]{li} there is a natural isomorphism
$A\otimes_B M\cong\projlim{i} M/\m_B^i.$ Thus the first statement
follows from lemma \ref{m-adic}. The isomorphism in (ii) is a
standard fact of homological algebra using the flatness of $A$
over $B.$ If $M$ is finitely generated over $R$ we can find
distinguished polynomials $F_1,\ldots,F_n$ and a surjection
\[\bigoplus_{i=1}^rA/AF_i \twoheadrightarrow M.\] Denoting the
kernel of it by $N$ the long exact $\E^{\bullet}$-sequence (for
both $A$ and $B$) gives
\begin{eqnarray*}
0=\E^0(N)\to\E^1(M)\to \bigoplus_{i=1}^r\E^1(A/AF_i) \to
\E^1(N)\to\E^2(M)\to 0
\end{eqnarray*}
\mbox{ and} $\E^i(N)\cong\E^{i+1}(M) \mbox{ for all } i\geq 3.$ An
easy calculation using the ``right" version of the Weierstrass
preparation theorem shows that
\[\E^1_A(A/AF_i)\cong A/F_iA\cong B/F_iB\cong \E^1_B(B/BF_i)\] as right modules. Thus
$E^1(M)$ is finitely generated over $R$ and $E^1_B(M)$ bears a
natural right $A$-module structure by lemma \ref{A-sub} (indeed
its right version). Applying this also to $N$ instead
of $M$ we obtain the statement (ii) by induction. Now (iii) follows from (i) and (ii). 
\end{proof}

 Assume now that
 \begin{itemize}
 \item $R$ is still a compact noetherian local ring with maximal ideal $\m_R$ as before, but
 now such that
 \item  $gr_{\m_R}R$ is a commutative noetherian regular equidimensional integral
 domain which is an algebra over some field and
 \item  $\sigma$ induces the identity on the residue class field $k,$ i.e.
 $gr_{\m_A}A=gr_{\m_B}B$ is also  commutative.
 \end{itemize}
   Then, by \ref{grR} and its
proof,  both $A$ and $B$ are Auslander regular and without zero
divisor. Thus there is a dimension theory for $A$- and
$B$-modules, for a thorough treatment of which we refer the reader
to \cite{co-sch-su} or \cite{ven1}. The previous lemma shows that
the canonical (co-)dimension filtrations of $M\in
\mbox{$A$-mod}^R$ considered as $A$- or $B$-module coincide.
Recall that  $M$ is called {\em pseudo-null} if $j_A(M)\geq 2,$
i.e.\ if $\E_A^0M=\E_A^1M=0$  (similarly for $B$).  Since $R$ is
without zero divisors the set $S:=R\setminus\{0\}$ is an Ore set
of $R$ and thus the total ring of quotients $Q:=R_S$ exists. By
\cite[prop. 2.12]{cohn} $\alpha$ and $\delta$ extend uniquely to
$Q,$ thus  $S$ is also an Ore Set for $B$ and there is a natural
isomorphism of rings
\[C:=B_S=Q[X;\sigma,\delta]\] (see below for further explanation in the case where
$\delta=\sigma -\mathrm{id}$). Then it is well-known that $C$ is a
principal ideal domain (division algorithm). Now we are able to
prove the following characterization of pseudo-null $A$-modules
which are finitely generated over $R:$

\begin{prop}\label{pseudo-torsionII} Let $M$ be in $\mbox{$A$-mod}^R.$ Then the following holds:
  $M$ is pseudo-null if and only if $\rk_R M=0,$ i.e.\ if and only if $M$ is a torsion
$R$-module.
\end{prop}

\begin{proof}
Let $M$ be a pseudo-null $A$-module which is finitely generated as
$R$-module. Then
\[\E^i_C(C\otimes_B M)\cong\E^i_B(M)\otimes_B C=0\] for all $i$ because $C$ has global
dimension less or equal to $1.$ It follows that $Q\otimes_RM\cong
C\otimes_BM=0,$ i.e.\ $M$ is $R$-torsion. Since $gr_{\m_A}A$ is a
commutative noetherian regular equidimensional integral  domain
which is an algebra over some field, the other direction follows
as in \cite[lem. 6.4]{ochi-ven} using Hilbert-Samuel polynomials.
\end{proof}

We write   $\mathcal{C}$ for the full subcategory of $A$-mod
consisting of pseudo-null modules. Since $\mathcal{C}$ is a Serre
subcategory the (abelian) quotient category $\mbox{$A$-mod}/
\mathcal{C}$ exists together with a natural functor
\[q:\mbox{$A$-mod}\to \mbox{$A$-mod}/ \mathcal{C}.\] By a sort of
homomorphism theorem for abelian categories one may identify the
analogous quotient category of  $\mbox{$A$-mod}^R$ with a full
subcategory of the prior one \[ \mbox{$A$-mod}^R/\mathcal{C}\cap
\mbox{$A$-mod}^R\subseteq \mbox{$A$-mod}/ \mathcal{C}.\]
Henceforth we will denote $ \mathcal{C}\cap \mbox{$A$-mod}^R$ just
by $\mathcal{C}^R.$  On the other hand restriction of scalars and
tensoring with $C\otimes_B -$ gives rise to an exact functor of
abelian categories
\[\Phi : \mbox{$A$-mod}^R\to \mbox{$B$-mod}^R\to \mbox{$C$-mod},\] the kernel of which is
just $\mathcal{C}^R$ by the above proposition. By the universal
property of  quotient categories, $\mbox{$A$-mod}^R/
\mathcal{C}^R$ can be identified with a full subcategory of
$C$-mod, which we call by abuse of notation $ \mbox{$C$-mod}^R.$
In terms of left ideals,i.e.\ of cyclic modules this can be
rephrased as follows: Let $\Ref$ denote the set of non-zero left
ideals $I$ such that $M:=A/I$ is in $\mbox{$A$-mod}^R$ (i.e.\ $I$
contains a distinguished polynomial) and  does not contain any
non-trivial pseudo-null submodule. By \cite[lem.\ 3.12]{co-sch-su}
$\Ref$ coincides with the set  of all non-zero reflexive left
ideals $I$ of $A$ such that $I$ contains a distinguished
polynomial. For all $I$ in $\Ref,$ $I_B:=I\cap B$ is a  non-zero
left ideal of $B$ and its legalization $I_C:=(I_B)_S=CI_B$ is  a
non-zero left ideal of $C.$ Writing  $\W$ for the set of all
non-zero left ideals of $C$ which contain a distinguished
polynomial (with coefficients in $R$)   we get a map
\begin{eqnarray*}
(-)_C:\Ref&\to& \W \\ I&\mapsto & I_C.
\end{eqnarray*}

There is also a map in the reverse direction $(-)_A: \W\to\Ref,$
induced by
\[J\mapsto J_B:=J\cap B\mapsto J_A:=AJ_B.\]
It is easy to see that  $(J_A)_C=J$ for all $J\in\W.$ Indeed,
since the left ideals $J_B\subseteq J_A$ contain a  WP the
$B$-module $B/J_B$ is the restriction of the $A$-module $A/J_A$ by
lemmas \ref{A-sub} and  \ref{pseudo-torsionI} (i), i.e.\
$J_B=J_A\cap B=(J_A)_B.$ By \cite[prop. 1.15]{mc-rob} it follows
that $J=CJ_B=(J_A)_C.$ On the other hand, it holds $I=(I_C)_A$ for
any $I$ in $\Ref:$ using (loc.\ cit.) again one obtains an
inclusion $I_B\subseteq CI_B\cap B=(I_C)_B$ with becomes an
equality after scalar extension to $C.$ Thus $(I_C)_A/I\cong
(I_C)_B/I_B$ is pseudo-null by prop.\ \ref{pseudo-torsionII} and
hence zero by assumption. Thus we have proved the following

\begin{lem}
 $(-)_C$ and $(-)_A$ are inverse maps, i.e.\ they are bijections
 of $\Ref$ and $\W.$
\end{lem}

We shall write $\W_{m}$ and $\W_{d}$ for the subsets of $\W$
consisting of those non-zero left ideals of $C$ which are
generated by a monic or  distinguished polynomial  in $B,$
respectively. I am very grateful to R. Sujatha and J. Coates for
pointing out to me the following consequence of the above lemma
and the resulting question below.

\begin{cor}
An left ideal $I\vartriangleleft_l  A$ is principal if and only if
$I_C$ belongs to $\W_{m}.$ Furthermore, the sets $\W_{m}$ and
$\W_{d}$ coincide.
\end{cor}

\begin{proof}
First we observe that for a monic polynomial $f$  in $B$ the
following holds: \[Cf\cap B=Bf.\] Indeed, it is easy to see that
for monic polynomials there is an Euclidean algorithm both in $B$
and $C, $ which implies the claim by the flatness of $Q$ over $R.$
Hence, if $J$ is in $\W_{m},$ say $J=Cf,$ with $f\in B$ monic, we
obtain $J_B=Bf,$ and so $J_A=Af=Ah,$ where $h$ is a distinguished
polynomial, by the Weierstrass Preparation Theorem. Hence
$J=(J_A)_C=Ch,$ i.e.\ $J$ belongs to $\W_{d}.$ If we take $J=I_C,$
it follows that $I=J_A=Ah$ is principal. The other implications of
the corollary are obvious.
\end{proof}

Hence we are left with the following key question concerning the
principality of ideals in $\Ref$

\begin{ques}
Is every non-zero left ideal of $C$ which contains a distinguished
polynomial generated by a monic polynomial in $B?$
\end{ques}

This seems to be a basic question. The answer is yes if and only
if every ideal $I$ in $\Ref$ is principal.

Now, we return to the above quotient category: since simple
objects in $\mbox{$A$-mod}^R/ \mathcal{C}^R$ can be represented by
some cyclic module $A/I,$ one immediately obtains

\begin{thm}\label{maxleft-simple}
There exists an equivalence of categories \[\mbox{$A$-mod}^R/
\mathcal{C}^R \cong \mbox{$C$-mod}^R,\] under which the
equivalence classes of simple objects of the quotient category
correspond uniquely to equivalence classes of maximal left ideals
of $C$ which contain a distinguished polynomial. Furthermore, $
\mbox{$C$-mod}^R$ is the smallest abelian subcategory of $C$-mod
which is closed under extensions and contains the simple modules
corresponding to the latter class of maximal left ideals.
\end{thm}

Here the equivalence classes in the quotient category are just the
classes of to each other isomorphic objects, while two left ideals
$J_1,\; J_2$ of $C$ are equivalent if and only if there exists an
isomorphism $C/J_1\cong C/J_2$ of $C$-modules.

We end this section by applying the above result to prime ideals.
For simplicity we assume now that $\delta=\sigma-\text{id}.$ For
any subring $\tilde{B}\subseteq A$ containing the ring
$B:=R[Y;\sigma,\delta]\cong R[Z;\sigma]$ ($Y=Z+1$) we write
$Spec_R(\tilde{B})$ for the set of those prime ideals $\p$ of
$\tilde{B}$ such that $\tilde{B}/\p$ is a finitely generated
$R$-module. Note that the zero ideal is {\em not} contained in
$Spec_R(\tilde{B}).$

\begin{lem}
Let $\p$ be in $Spec_R(A).$ The intersection with $\tilde{B},$
$\p\mapsto \p\cap \tilde{B}$ induces an injective map
\[Spec_R(A)\to Spec_R(\tilde{B}).\]
\end{lem}

\begin{proof}
For $\p\in Spec_R(A),$ $\tilde{B}/\p\cap \tilde{B}=A/\p$ by the
Weierstrass preparation theorem. Since $\tilde{B}$ is dense in $A$
there is only one $A$-module structure on $\tilde{B}/\p\cap
\tilde{B}$ extending the $B$-module structure which implies
injectivity.
\end{proof}

Now we specify $\tilde{B}$ to be the skew Laurent polynomial ring
\[\tilde{B}:=R[Z,Z^{-1};\sigma].\] Since the ring-automorphism extends uniquely to the skewfield of
fractions $Q$ of $R,$ the ring $\tilde{B}$ embeds into
\[\tilde{C}:=Q[Z,Z^{-1};\sigma]\] which is - as legalization of
$Q[Z,\sigma]$ - a principal ideal domain. Moreover, \linebreak
$Q[Z,Z^{-1};\sigma]$ is also the quotient ring  $\tilde{B}_S$ of
$\tilde{B}$ with respect to the set $S=R\setminus \{0\}$ of
regular elements of $R.$ Indeed, since $\sigma$ is an
automorphism, the elements of $R[Z,Z^{-1};\sigma]$ and
$Q[Z,Z^{-1};\sigma]$ can be written either with coefficients on
the left or right site of the powers of $Z.$ Let $\p$ be in
$Spec_R(\tilde{B})$ and  assume that $\tilde{B}/\p$ has no
$R$-torsion. Then, $\p\cap R=0$ and thus $\p_S\subsetneqq
\tilde{B}_S= C$ is a proper prime ideal, see \cite[prop.
2.1.16]{mc-rob}. Denoting by $Spec_{R,tf}(A)$ the subset of
$Spec_R(A)$ consisting of those prime ideals $\p$ such that $A/\p$
is without $R$-torsion we obtain

\begin{prop}\label{spec-tf}
Under the above conditions there is a natural embedding
\[Spec_{R,tf}(A)\subseteq Spec_{R,tf}(\tilde{B})\subseteq
Spec(\tilde{C})\setminus\{0\}(\subseteq
Spec(Q[Z;\sigma])\setminus\{0\}).\]
\end{prop}

By $G(A)$ and $G_R(A)$ we denote the group of (fractional)
$c$-ideals and the subgroup which is generated by $c$-prime ideals
in  $Spec_{R,tf}(A),$ respectively. For the definition and a
review of its basic properties we refer the reader to
\cite[{\S}4]{co-sch-su}. Note that under our current assumptions $A$
is a maximal order (loc.\ cit.\ lem. 2.6). Since for any maximal
order $\mathcal{O},$ $G(\mathcal{O})$ is the free abelian group on
the set of all prime $c$-ideals,  the proposition implies

\begin{cor}
There is an injective group homomorphism
\[G_R(A)\hookrightarrow G(Q[Z,Z^{-1};\sigma]).\]
\end{cor}

\section{Completely faithful objects}\label{com-faithful}

In this section we are going to apply the results of the previous
ones to certain Iwasawa algebras. In particular we are interested
in the existence of bounded and faithful objects in the quotient
category.

First we recall some general definitions. For this purpose let $A$
be any maximal order satisfying that every finitely generated
torsion $A$-module induces an object of finite length in the
quotient category $\mbox{$A$-mod}/ \mathcal{C}.$ The annihilator
ideal of an object $\mathcal{M}$ of  $\mbox{$A$-mod}/ \mathcal{C}$
is defined by $Ann(\mathcal{M}):=\sum_{q(M)\cong\mathcal{M}}
Ann_A(M),$ where $q$ denotes again the natural functor
\[q:\mbox{$A$-mod}\to \mbox{$A$-mod}/\mathcal{C}.\] Note that by
\cite[lem 2.5]{rob} in conjunction with \cite[cor of thm.
2.5]{chamarie2} $Ann (q(M))=Ann_A (M/M_{ps})$ where $M_{ps}$
denotes the maximal pseudo-null submodule of $M.$ Recall that the
object $\mathcal{M}$ is called {\em completely faithful} if
$Ann(\mathcal{N})=0$ for any non-zero subquotient $\mathcal{N}$ of
$\mathcal{M}.$ It is called {\em bounded} if $Ann(
\mathcal{N})\neq 0$ for any subobject
$\mathcal{N}\subseteq\mathcal{M}$ (In particular, the image $q(M)$
of any pseudo-null module $M$ is both bounded and completely
faithful, by definition).

We should also mention that  for a $A$-torsion module $M,$  its
image $\mathcal{M}:=q(M)$ decomposes uniquely into a direct sum
$\mathcal{M}\cong\mathcal{M}_{cf}\oplus\mathcal{M}_b$ where $
\mathcal{M}_{cf}$ is completely faithful and $\mathcal{M}_b$
bounded, cf. \cite[prop. 4.1]{co-sch-su}. Moreover,
$\mathcal{M}_{cf}$ or more general all completely faithful objects
of finite length $\mathcal{M}$ of $q( \text{$A$-mod})$ are cyclic,
i.e.\ $\mathcal{M}\cong q(\La/L)$ for some non-zero left ideal
$L\subseteq \La$ (loc.\ cit.\ lemma 2.7 and prop. 4.1).

Now let $G=H\rtimes \Gamma$ be an uniform pro-$p$-group as in
example \ref{semi-gen} and set $R:=\La(H)$ as well as
\[A:=\La(G)\cong \La(H)[[Y;\sigma,\delta]]=
R[[Y;\sigma,\delta]].\] For some left ideal $I$ of $A,$ consider
the  cyclic $A$-module $M=A/I$ and assume that it is finitely
generated as $R$-module and has no non-zero pseudo-null submodule.
Then, by proposition \ref{pseudo-torsionII}, $M$ is  torsionfree
as $R$-module. We shall  denote the annihilator ideal of $M$ by
$J=Ann_A(M).$

\begin{lem}
In the above situation, assume that $M$ is bounded, i.e.\ $J\neq
0.$ Then the module $\Lambda/J$ is finitely generated over $R,$
too, and has no non-zero pseudo-null submodule, either.
\end{lem}

\begin{proof}
We use the same notations as in lemma \ref{ann-intersection} and
\ref{intersection} ($\Lambda$ being replaced by $A$). Since
$I_n=Ih_n^{-1}$ and $h_n$ is a unit in $A$, obviously $A/I_n \cong
A/I$ and we obtain an injection

\[A/J_i\subseteq\prod_{n=0}^i A/I_n\cong \prod M.\]
Thus, for all $i\geq 0,$ the module $A/J_i$ has no non-zero
pseudo-null submodule. Now, assume that infinitely many
subquotients $J_i/J_{i+1}$ of the chain of left ideals

\[\cdots\subseteq J_{i+1}\subseteq J_i\subseteq \cdots \subseteq J_1\subseteq J_0=I\] are non
trivial (here we assume without loss of generality that $h_0=1,$
i.e.\ $J_0=I_0=I$). Then, since $0\neq J_i/J_{i+1}\subseteq
A/J_{i+1}$ cannot be pseudo-null, the image $q(J_i/J_{i+1})$ would
be non-zero  for infinitely many $i$ contradicting the fact that
all torsion $A$-modules have finite length in the quotient
category. Thus the above chain becomes stationary, say $J=J_k,$
and hence \[A/J=A/J_k\subseteq \prod_{n=0}^k M\] is finitely
generated over $R$ and without non-zero pseudo-null submodule.
\end{proof}

This lemma tells us  that the identification

\[\mbox{$A$-mod}^R/ \mathcal{C}^R \cong
\mbox{$C$-mod}^R\]

in theorem \ref{maxleft-simple} respects ``boundedness", i.e.\ an
object $\mathcal{M}=q(M)$ is bounded if and only if $C\otimes_B M
$ is a bounded $C$-module. In particular, $q(A/I)$ is bounded if
and only if $I_C$ is a bounded ideal of $C,$ i.e.\ contains a
nonzero ideal. Also, using \cite[prop. 4.4]{co-sch-su}, it implies
that $G_R(A)\cong G((\mbox{$A$-mod}^R/\mathcal{C}^R)^b)$ where the
latter denotes the Grothendieck group of the full subcategory of
bounded objects  $(\mbox{$A$-mod}^R/\mathcal{C}^R)^b$ of $
\mbox{$A$-mod}^R/\mathcal{C^R}.$ Hence, by the results of the
previous section we obtain the following

\begin{prop}
There is an injective group homomorphism
\[ G((\mbox{$A$-mod}^R/\mathcal{C}^R)^b)\hookrightarrow G(Q[Z,Z^{-1};\sigma]) .\]
\end{prop}

There are two special cases of particular interest. First, if $G$
is the semidirect product of two copies of $\zp,$ we have seen
that $ Q[Z,Z^{-1};\sigma]$ is simple (cf. remark \ref{rem-ann}),
thus $G(Q[Z,Z^{-1};\sigma])=0$ and we obtain the  following

\begin{thm}\label{th1} Let $G$ be a non-abelian group isomorphic to $\zp\rtimes\zp.$
Then the image $q(M)$ of any $M\in \text{$A$-mod}^R,$ is
completely faithful, i.e.\ the subcategory $q( \text{$A$-mod}^R)$
of $\mbox{$A$-mod}/\mathcal{C}$ consists only of completely
faithful objects. In particular,  all objects $\mathcal{M}$ of
$\phantom{,} q( \text{$A$-mod}^R)$ are cyclic, i.e.\
$\mathcal{M}\cong q(\La/L)$ for some non-zero left ideal
$L\subseteq \La.$
\end{thm}


Of course, the theorem is also an immediate consequence of
corollary  \ref{cor-faithful} and proposition \ref{dim}.

 Second, if $G$ is an open subgroup of $GL_2(\zp)$ such that $G=H\times\Gamma$
where $H=G\cap SL_2(\zp)$ and $\Gamma=\text{center}(G),$ then
$\sigma$ is the identity, i.e.\ the variable $Z$ or $Y$ commute
with $Q.$ It is well known and easy to check that then every
two-sided   ideal of $Q[Y]$ can be generated by an element in
$Z(Q)[Y]\subseteq Q[Y]$ where $Z(Q)$ denotes the center of the
skewfield $Q.$ In  \cite{howson2001} S. Howson has proved that
$Z(R)=\zp$ in this case. In order to determine the bounded
$A$-modules which are finitely generated over $R$ and not
pseudo-null it might be crucial to answer the following

\begin{ques}
Let $H$ be any open pro-$p$-subgroup of $SL_2(\zp).$ Is the center
$Z(Q)$ of the skewfield of fractions $Q$ of $\Lambda(H)$ (i) equal
to $\qp,$ (ii) algebraic over $\qp$ or (iii) has it transcendent
elements over $\qp?$ Are the elements of $Q\setminus Z(Q)$
algebraic or transcendent over $Z(Q)?$
\end{ques}

 In any case, we know that the simple objects which can be represented by  $A$-modules which are
  finitely generated over $R$ correspond to those maximal left ideals of $Q[Y]$ which contain a
  distinguished polynomial. If item (i) of the question holds it would probably tell us that
   the only bounded   objects are those which are induced from the center of $G.$

 Since $\La(H)\subseteq\La(G)$
for the canonical normal subgroup $H$ of $G,$ we may consider any
$\La(G)$-module $M$ also as $\La(H)$-module. If such $M$ is even
finitely generated over $\La(H),$ then it is known that
$\dim_{\La(H)}M\geq\dim_{\La(G)}$ (\cite[lem. 6.4]{ochi-ven}). We
now will prove that in the case where $G\cong\zp\rtimes\zp.$
equality holds.

\begin{prop}\label{dim} Let $G\cong \zp\rtimes\zp.$
With the above notations it holds that
$\dim_{\La(H)}M=\dim_{\La(G)}.$ In particular, $M$  is a
pseudo-null \La(G)-module if and only if it is a torsion
$\La(H)$-module.
\end{prop}

 Thus a pseudo-null $\La(G)$-module which is finitely generated over $\La(H)$ is - as $\La(H$)-module up to
 finite modules  - the direct sum of a finitely generated \zp-module and summands of the form
 $\mathbb{Z}/p^n[[H]].$

\begin{proof}
Since $M$ is zero-dimensional (as \La(G)- or \La(H)-module) if and
only if it is finite, we only have to prove

{\em Claim:} If $M$ is torsion-free as \La(H)-module than it is
not a pseudo-null module.

This is just the statement of proposition \ref{pseudo-torsionII},
but we want to give a different argument which might be
interesting on its own: Since $X$ is a normal element of \La(G),
$M$ has a descending chain of submodules
$\{X^nM\}_{n\in\mathbb{N}}.$ If $M$ is torsionfree as
$\La(H)\cong\zp[[X]]$-module , then it es easy to see that the
subquotients $X^{n+1}M/X^nM$ are isomorphic to $M/XM$ as
\zp-module. Since the latter module cannot be finite by the strong
Nakayama lemma, it follows that these subquotients have dimension
$1$ both over \La(G)\ and  \La(H)\ (see \cite[prop. 3.5 (ii),cor.
4.8]{ven1}). Choosing a good filtration on $M$ (and taking the
induced filtrations on $X^nM$) we get an infinite descending chain
of associated graded modules
\[grM\supsetneqq gr(XM)\supsetneqq \ldots \supsetneqq gr(X^nM) \supsetneqq \ldots \] with
subquotients $gr(X^nM)/gr(X^{n+1}M)\cong gr(X^nM/X^{n+1}M)$ of
Krull dimension $1$ over $gr\La(G)\cong\mathbb{F}_p[X_0,X_1,X_2]$
by \cite[thm. 3.22]{ven1}. This implies that $grM$ has Krull
dimension $2,$ i.e. $\dim_{\La(G)}M=2,$ by the next lemma.
\end{proof}

\begin{lem}
Let $R$ be a commutative Noetherian ring with finite Krull
dimension, $M$ a finitely generated $R$-module which has an
infinite chain of submodules \[ M=M_0\supsetneqq M_1\supsetneqq
\ldots \supsetneqq M_i \supsetneqq \ldots\] with non-zero
subquotients $M_i/M_{i+1}$ of Krull-dimension $\alpha.$  Then, the
Krull dimension of $M$ is (strictly) greater than $\alpha.$
\end{lem}

\begin{proof}
Assume the contrary. Then $deg(M):=\sum_{\dim
R/\p=\alpha}length_{A_\p} M_\p$ is finite (compare \cite[VIII {\S}
1.5]{bourbaki8}). Since the function $deg$ is additive (loc.cit.),
all subquotients contribute a nonzero positive number to $deg(M)$
by assumption, a contradiction.


\end{proof}

\begin{rem}
The lemma corresponds to the more or less well known fact that for
a commutative Noetherian ring the classical Krull-dimension
coincides with the Gabriel-dimension or, depending on the
definition, with the Gabriel-dimension minus one, see the notes
``Krull dimension" in \cite[p. 239]{GoWa}.
\end{rem}

\begin{ques}
For which $p$-adic Lie groups without $p$-torsion, does it hold
that the Gabriel-dimension of every finitely generated
$\Lambda(G)$-module coincides with its Bjoerk-dimension,
respectively the Krull-dimension of the associated graded module?
\end{ques}

Of course the question is answered positively for all {\em
abelian}  $p$-adic Lie groups without $p$-torsion.

\section{Unique Factorisation Rings}\label{ufd}

Chatters and Jordan studied in \cite{cha-jor} a class of
Noetherian (= left and right Noetherian) rings which  in the
commutative case consists of the Noetherian unique factorisation
domains:

{\em A prime Noetherian ring $A$ is called {\em unique
factorisation ring (UFR)} if every non-zero prime ideal of $A$
contains a non-zero principal prime ideal.}

If $A$ satisfies the descending chain condition on prime ideals
(e.g.\ if the classical Krull dimension $\dim (Spec(A))$ of $A$ is
finite, where $Spec(A)$ is the poset of prime ideals, see
\cite[Appendix 3.]{GoWa}), then the above condition is equivalent
to

{\em Every height-$1$ prime ideal of  $A$ is principal.}

The meaning of {\em principal ideal $I$} is that $I=aA=Aa$ for
some element $a \in A.$ Such elements are called {\em normal} and
if $I$ is a prime ideal in addition than $a$ is called {\em prime}
element. The class of UFRs contains the following much more
restricted  one (Remark (5) in \cite{cha-jor}):

{\em A Noetherian integral domain $A$ which contains at least one
height-$1$ prime ideal is called {\em unique factorisation domain
(UFD) } if every height-$1$ prime $P$ of $A$ is principal and
$A/P$  is an integral domain.}

We also should mention that any UFR $A$ is a unique factorisation
ring in the sense of \cite{akmu}, i.e.\ $A$ is a maximal order
and any prime {\em $c$-ideal}\ \ - a prime ideal which is (left or
right) reflexive - is principal (cf.\ the remark after prop. 1.5
in \cite{akmu}).

Now we come back to our example \ref{semi}, where $A=\fp \kl G
\kr$ is the completed group algebra of the semidirect product
$G=\zp\rtimes\zp.$ It seems that $A$ is the first example of a
complete respectively power series ring which admits unique
factorisation.

\begin{thm}\label{ufd-thm}
The ring $A=\fp \kl G \kr=\fp [[X,Y;\sigma,\delta]]$ is a UFD.
Furthermore, $(X):=AXA$ is the only prime ideal of height one in
$A$ and $X$ is - up to multiplication by units - the only prime
element of $A.$
\end{thm}

\begin{proof}
Since $(X)$ is the kernel of
\[ \bar{ }:\fp [[X,Y;\sigma,\delta]]\to \fp[[Y]]\]
it is a prime ideal which coincides with $XA$ and by symmetry also
with $AX$ (if the coefficients are written on the right of the
powers of $Y$). By the principal ideal theorem \cite[thm.
4.1.11]{mc-rob}  $(X)$ has height $1.$ Now assume that there is
some height-$1$ prime ideal $P$  of $A$ different from $(X).$ Then
$P$ contains an element $p$ with $\bar{p}\neq 0$ in $\fp[[Y]].$ By
the Weierstrass preparation theorem $p$ is in $R[Y;\sigma,\delta]$
up to a unit where $R=\fp[[X]].$ The remark \ref{rem-ann} (ii)
implies that $P\cap R\neq 0,$ i.e. $X^n\subseteq P\cap R\subseteq
P$ for some $n$ because all ideals of $R$ are of the form $(X^n).$
Since $P$ is prime it follows that $(X)\subseteq P,$ a
contradiction. Thus the unique height-$1$ prime $(X)$ is principal
with $A/(X)$ being integral, i.e.\ $A$ is a UFD.
\end{proof}

Similarly, it is seen that $X^n,\; n\in\mathbb{N},$ are - up to
multiplication by units - the only normal elements of $A.$

Let $C(X)$ be the set of elements of $A$ which are regular modulo
$(X)$ and $S=\{X^n| n\in\mathbb{N}\}.$ Then it is known (and in
this example easy to verify) that the classical (left and right)
localisations $A_{(X)}:=A_{C(X)}$ and $A_S$ of $A$ with respect to
$C(X)$ and $S$ exist. They have the following properties

\begin{enumerate}
\item $A_{(X)}$ is a bounded local principal ideal domain with a single prime ideal $XA_{(X)}$
and whose only (left or right or two-sided) ideals are of the form
$X^nA_{(X)}$ (see \cite{cha} prop. 2.8 and its proof in
conjunction with th. 2.7, compare also with \cite[prop. 2.5
]{chamarie1} ).
\item $A_S$ is simple (\cite[lem. 2.1]{cha-jor}).
\item $A=A_{(X)}\cap A_S$ (\cite[thm. 2.3]{cha-jor}). This should be compared also with
Chamarie's result 3.3-3.5 in \cite{chamarie2}: While $A_S$ is
"responsible" for the faithful modules the bounded ring $A_{(X)}$
covers the bounded modules (in the quotient category of $A$-mod
modulo the Serre subcategory of pseudo-null $A$-modules, see
\cite{co-sch-su} or \cite{ven1}).
\end{enumerate}

\begin{rem}
Let $A=\fp \kl G \kr$ and $J$ be a left ideal of $A$ such that
$M:=A/J$ does not contain any pseudo-null submodule. Then $M$ is
faithful in the quotient $\mbox{$A$-mod}/\mathcal{C}$ if and only
if no power $X^n,$ $n\in\mathbb{N},$ is contained in $J.$ Indeed,
under the assumptions the annihilator of the image of $M$ in the
quotient category coincides with the usual annihilator ideal
$Ann_A(M).$ If $M$ is bounded, the latter ideal is a $c$-ideal in
$A,$ i.e.\ a power of the unique prime ideal $(X)$ (\cite[lem
4.3]{co-sch-su}). In particular, $X^n$ is contained in $J$ for
some $n.$ On the other hand, assume that   $X^n\in J.$  Then $J$
contains obviously the two-sided ideal $AX^n.$
\end{rem}

We should also mention that the fact that $\fp \kl G \kr$ is a UFR
in the sense of \cite{akmu} can be derived more easily from the
following general theorem:

\begin{thm}
Let $A$ be a local noetherian integral ring which is a maximal
order (in its ring of quotients). If  the  global homological
dimension of $A$ is less or equal to $2$ then $A$ is a UFR in the
sense of \cite{akmu}.
\end{thm}

\begin{proof}
Conferring prop. 1.8. (2) in \cite{akmu} it is sufficient to prove
that every reflexive ideal $I$ is principal: Denote by
$I^+:=\Hom_A(I,A)$ the dual of $I$ and choose a projective
resolution of $I^+$

\[ 0\to P_1\to P_0 \to I^+\to 0.\]
Applying the functor $-^+$ we obtain  an exact sequence

\[0\to I\to P_0^+\to P_1^+\to \mbox{Ext}^1_A(I^+,A)\to 0\]

because $I\cong I^{++}$ is reflexive. Now it follows by
homological algebra that $I$ has projective dimension equal to
zero (since the projective dimension of $ \mbox{Ext}^1_A(I^+,A)$
is less or equal to $2$). But as $A$ is local it follows that $I$
is a free left $A$-module. By symmetry it is also a right
principal ideal and the theorem follows.
\end{proof}
Recall that a {\em principal ideal domain (PID)} is an integral
ring all left and right ideals of which are principal.

\begin{cor}
Let $R$ be a local PID and $\pi$ a generator of the unique maximal
ideal $\m=(\pi)$ of $R$ and assume that $R$ is complete in its
$\pi$-adic topology. If $\sigma$ is a ring endomorphism of $R$
which induces an automorphism $\bar{\sigma}$ of $\kappa:=R/(\pi),$
then, for any $\sigma$-derivation $\delta$ as in section
\ref{powerseries}, the skew power series ring
$R[[Y;\sigma,\delta]]$ is an UFD in the sense of  \cite{akmu}.
\end{cor}

\begin{proof}
 Since $\kappa$ is a skew field its global dimension is zero. Hence, by corollary
 \ref{grR} and \cite[thm. 7.5.3 (iii)]{mc-rob},
 $\mathrm{gl}A\leq \mathrm{gl}gr_\pi R+1=\mathrm{gl}\kappa[x]+1=2$ and the result follows from
 proposition \ref{local}, corollary \ref{grR} (iii)   and the  theorem.
\end{proof}

In general, there is a canonical homomorphism
\[G(R[Y;\sigma,\delta])\to G(R[[Y;\sigma,\delta]]),\] induced by
$I\to \widehat{I}:=A\otimes_B I\cong I\otimes_B A.$ In those cases
where this map is surjective (e.g.\  for $\fp\kl G\kr$ as above)
the UFR property for $B$ in the sense of \cite{akmu} implies the
UFR property for $A$ (supposed the latter ring is a maximal
order). Anyway, one obtains the following implication

\begin{prop}
Let $R$ be a compact noetherian local ring with maximal ideal
$\m_R.$  If $R[Y;\sigma,\delta]$ is a UFR in the sense of
\cite{akmu}, then all ideals in $G_R(A)$ are principal.
\end{prop}

We should mention that  theorem  2.4 and proposition 1.8 (2) of
\cite{akmu} imply the following: If a noetherian ring $R$ is a UFR
in their sense, then so is $R[Z,\sigma].$  E.g., if $\Lambda(H)$
is a UFR, then all ideals in $G_{\La(H)}(\La(G))$ are principal,
where $H\subseteq G$ are the groups of example \ref{semi-gen} or
in particular of example \ref{gl2} where $H\subseteq G\subseteq
GL_2(\zp). $

\begin{proof}
Let $I$ be a prime $c$-ideal of $A$ such that $A/I$ is finitely
generated over $R,$
 set $I_B:=B\cap I$ as before and consider the following exact and commutative diagram
   \[ \xymatrix{ 0\ar[r] &(I_B)^{++}\ar[r] &{B^{++}}\ar[r] &{\E^1\E^1(B/I_B)}  \\
              0\ar[r] &I_B\ar[r]\ar^{\varphi_{I_B}}[u] &{B}\ar[r]\ar@{=}^{\varphi_{B}}[u] &B/I_B\ar[r]\ar[u] &0 \\
                 }\]
 Here, the first line is easily obtained by forming twice the long exact
 $\E^\bullet$-sequence of the bottom line and the vertical homomorphisms are  (induced by)
 the canonical maps from a module to its bi-dual. By the snake lemma, lemma \ref{A-sub} and
 \ref{cor-faithful} (i), the functor $A\otimes_B-$ is faithfully flat on the cokernel $W$ of
 $\varphi_{I_B}.$ But  $A\otimes_B W\cong \mathrm{cokern}(\varphi_{I})$ (as $A$-module) is zero
 by assumption. Thus, by left-right symmetry, $I_B$ is also a reflexive ideal and maps onto
 $I,$ i.e.\  the following restriction of the above homomorphism is surjective
\[G_R(B)\twoheadrightarrow G_R(A).\] The result follows.
\end{proof}

In section \ref{faithfulmod} we have seen that $\La(G)$ possesses
an abundance of faithful modules. Now we will describe an infinite
series of bounded modules which induces an infinite series of
pairwise non-isomorphic simple objects in the quotient category by
pseudo-null modules.  Indeed, by \cite[prop. 4.4]{co-sch-su} the
prime $c$-ideals, i.e.\ those prime ideals which are reflexive
(both as left and right $\La(G$-)module), correspond to the
annihilator ideals of  bounded simple objects. Thus it is
sufficient to find an infinite series of prime ideals which are
principal and thus reflexive.

To this end we set $\omega_{-1}:=1,$  $\omega_n:=(X+1)^{p^n}-1,$
$n\geq 0,$ and denote the cyclotomic polynomial by
\[\xi_n:=\frac{\omega_n}{\omega_{n-1}},\;\; n\geq 0.\] Then, $\xi_0=\omega_0=X,$ while
\[  \xi_n=\sum_{i=0}^{p-1} (1+X)^{ip^{k-1}},\;\; k\geq 1. \]

\begin{prop}
The polynomials $\xi_n$ are prime elements of $\Lambda=\La(G),$
i.e.\ $(\xi_n):=\Lambda\xi_n=\xi_n\Lambda$ is a prime ideal of
\La\ for every $n\geq 0$. Furthermore, the quotient
$\Lambda/(\xi_n)$ is isomorphic to the integral ring
$\zp[\zeta_{p^n}][[Y;\overline{\sigma},\overline{\delta}]],$ where
$\zeta_{p^n}$ denotes a primitive $p^n$-root of unit and
$\overline{\sigma}$ and $\overline{\delta}$ are induced by
$\sigma$ and $\delta,$  respectively.
\end{prop}

In fact, using the integrality of
$\zp[\zeta_{p^n}][[Y;\overline{\sigma},\overline{\delta}]]$ and
the change of rings spectral sequence for Ext-groups, one can
prove that the $\La$-modules
\[\Lambda/(\xi_n)\cong\zp[\zeta_{p^n}][[Y;\overline{\sigma},\overline{\delta}]]\] {\em themselves} induce
simple objects $q(\Lambda/(\xi_n))$ in the quotient category. Note
that, in general, a prime $c$-ideal $\p$ is only the largest
two-sided ideal in some reflexive left ideal $I$ such that
$q(\La/I)$ is simple (with annihilator ideal $\p$).

\begin{proof}
First, we show that $\omega_m$ and $\xi_n$ are normal elements in
$\La:$ Since $\omega_n\La=\prod \zp[[X]]\omega_n Y^i$ is the
kernel of the canonical map
\[{\mathbb{Z}_p}[[X,Y;\sigma,\delta]]\cong\La(G)\to\La((H/H_n)\rtimes \Gamma)\cong(\zp[[X]]/\zp[[X]]\omega_n)[[Y;\overline{\sigma},\overline{\delta}]],\]
this is a two-sided ideal, which by symmetry coincides with
$\La\omega_n,$ thus $\omega_n$ is normal. Hence, we
obtain\[\omega_{n-1}\xi_n\La=\omega_n\La=\La\omega_n=\La\omega_{n-1}\xi_n=\omega_{n-1}\La\xi_n,\]
i.e.\ $\xi_\La=\La\xi_n$ because \La\ does not have zero divisors.

Since $ (\zp[[X]]/\zp[[X]]\xi_n)[[Y;\sigma,\delta]]\cong
\zp[\zeta_{p^n}][[Y;\sigma,\delta]]$ is integral by cor. \ref{grR}
(note that $gr\zp[\zeta_{p^n}]\cong\mathbb{F}_p[\pi],$ where $\pi$
is the image of $X$ respectively $\zeta_{p^n}-1$).
\end{proof}

\section{The ``False Tate Curve"}\label{tate}

In this section we will study a Galois extension of a number field
$k$ which arises from a $p$-adic Galois representation which is
analogous to the local presentation induced by a Tate elliptic
curve. Thus, following J. Coates, we call it ``false Tate  curve."

So let $\alpha\in k^\times$ be a unit of $k$ which is not a root
of unity. Consider the discrete $G_k$-module
$A:=\bar{k}^\times/\Delta_\alpha,$ where $\bar{k}$ denotes a fixed
algebraic closure of $k,$ $G_k=G(\bar{k}/k)$ the absolute Galois
group of $k$ and $\Delta_\alpha=\alpha^\mathbb{Z}$ the subgroup of
$\bar{k}^\times$ which is generated by $\alpha.$ Then the
$p^n$-division points $A_{p^n}$ of $A$ contain the $p^n$-roots of
unity $\mu_{p^n}\subseteq \bar{k}^\times$ and we have natural
exact sequences of $G_k$-modules for every $n$
\[0\to\mu_{p^n}\to A_{p^n} \to \mathbb{Z}/p^n\mathbb{Z} \to 0,\]
where  the latter module has trivial $G_k$ action and where the
surjection $A_{p^n}\to\mathbb{Z}/p^n\mathbb{Z}$ is given as
follows: If $\beta\in \bar{k}^\times$ represents an element of
$A_{p^n},$ i.e.\ $\beta^{p^n}=\alpha^c$ for some $c\in\mathbb{Z},$
then its class $\beta\Delta_\alpha$ is send to $c\mbox{ mod }
p^n.$

Writing $T_p(B):=\projlim{n} B_{p^n}$ the Tate-module and
$V_p(B):=T_p(B)\otimes_\zp\qp$ for the $p$-adic representation
associated with an arbitrary discrete $G_k$ module $B,$ we obtain
the exact sequences
\[0\to \zp(1) \to T_p(A)\to \zp\to 0,\]
 where $\zp(1)=T_p(\mu)$ is the first Tate-twist of $\zp,$ and

\[0\to \qp(1)\to V_p(A) \to \qp\to 0.\]

If $G$ denotes the kernel of the action $\rho:G_k\to GL(V)$ of
$G_k$ on $V=V_p(A):$
\[G:=\ker(\rho:G_k\to GL(V));\]
then the field $k_\infty:=\bar{k}^G$ which is determined by $G$
via Galois theory is of course $k(A_{p^\infty})=k(\mu_{p^\infty},
\alpha^{p^{-\infty}}),$ the field, which arises by adjoining the
$p$-power roots of unity as well as $p$-power roots of $\alpha$ to
$k.$

For the rest of the paper, we assume that $k$ contains the $p$-th
roots of unity $\mu_p\subseteq k.$ By $k_{cycl}$ we denote the
cyclotomic $\zp$-extension $k_{cycl}=k(\mu_{p^\infty})$ of $k,$
by $H:=G(k_\infty/k_{cycl})$ and $\Gamma:=G(k_{cycl}/k)$ the
Galois groups of the extensions $k_\infty/k_{cycl}$  and
$k_{cycl}/k$ respectively. Then Kummer theory tells us that
\[\Hom(H,\mu_{p^\infty})\cong
\mbox{im}(\Delta_\alpha\otimes_\mathbb{Z}\qp/\zp\to
k_{cycl}^\times\otimes_\mathbb{Z}\qp/\zp)\cong \qp/\zp\] as
$G_k$-modules. In particular, $H\cong\zp(1)$ as $\Gamma$-module,
i.e.\ $G$ is isomorphic to the (non-abelian) semi-direct product
\[G\cong H \rtimes \Gamma =\zp(1)\rtimes \zp.\]

Of course, the arithmetic in $k_\infty/k$ depends decisively on
the choice of $\alpha.$ If we take $\alpha=p$ for example we  have
the following properties.

 \begin{lem}
 Let $k=\mathbb{Q}(\mu_p)$ and $k_\infty=k(\mu_{p^\infty},p^{p^{-\infty}}).$ Then the
 extension $k_\infty/k$ is \begin{enumerate}
 \item totally ramified at the unique place over $p,$ in particular there is just one prime
 of $k_\infty$ above $p.$
 \item  unramified outside $p.$
 \end{enumerate}
 \end{lem}

 \begin{proof}(The proof arose from a discussion with Y. Hachimori)
 It suffices to prove these statements for the finite extensions $k(\mu_{p^n},p^{p^{-n}}.)$
 Then item (ii) follows from \cite[Lemma 5]{birch}.

 Now consider the local extensions $K=\qp(\mu_{p^n})$ and $L=K(p^{p^{-n}})$ of $\qp.$ Since
 the extension $\qp(p^{p^{-1}})/\qp$ is not Galois, no $p$th root of $p$ \  can be contained in
 the cyclic extension $K/\qp.$ Hence, it follows from Kummer theory that the degree of $L$
 over $K$ is $[L:K]=p^n,$ i.e.\
 $[L:\qp]=[\mathbb{Q}(\mu_{p^n},p^{p^{-n}}):\mathbb{Q}](=(p-1)p^{2n-1})$ and in particular $p$
 does not split in $k(\mu_{p^n},p^{p^{-n}}).$  Since the maximal abelian quotient $G^{ab}$
  of $G=G(L/\qp)\cong
 G(L/K)\rtimes G(K/\qp)$ is isomorphic to \[G^{ab}\cong G(L/K)_{G(K/\qp)}\oplus
 G(K/\qp)=G(K/\qp)\]
 (note that $G(L/K)\cong \mathbb{Z}/p^n(1)$ has no non-zero $G(K/\qp)$-invariant quotient
because $G(K/\qp)$ acts via the cyclotomic character on $G(L/K)$),
the only cyclic extensions of $\qp$ in $L$ are contained in $K$
and cannot be unramified. Hence $p$ is totally ramified in
$k(\mu_{p^n},p^{p^{-n}})$ for all $n.$
 \end{proof}

Let $L$ be the maximal unramified abelian $p$-extension of
$k_\infty.$ Then its Galois group $X_{nr}=G(L/k_\infty)$ is a
compact $\Lambda(G)$-module where $\Lambda(G)$ denotes the Iwasawa
algebra \[\Lambda (G)=\zp\kl G \kr\cong
\zp[[X,Y;\sigma,\delta]].\] If we define $H_n$ and $\Gamma _n$ to
be the unique subgroup of  $H$ and $\Gamma$ of index $p^n$,
respectively, then $G_n:=H_n\rtimes\Gamma_n$ is an open normal
subgroup of $G$ of index $p^{2n}.$ Let us write
$k_{n,m}:=k_\infty^{H_n\rtimes\Gamma_n}=k(\mu_{p^{m+1}},\alpha^{p^{-n}}),$
 $0\leq n,m\leq \infty,$  for the intermediate fields of $k_\infty/k.$ Analogous to the module
 $X_{nr}$ we have the $\Lambda(\Gamma)$-modules $X_{nr}(k_{n,\infty}$ for every field
 $k_{n,\infty}=k_{cycl}(\alpha^{p^{-n}}).$ By class field theory there are canonical isomorphisms
 of $\Lambda(\Gamma)$- and $\Lambda(G)$-modules respectively:
 \begin{eqnarray*}
 X_{nr}(k_{n,\infty})&\cong&\projlim{ m} Cl(k_{n,m})(p)\\
 X_{nr}&\cong&\projlim{ n,m} Cl(k_{n,m})(p)\\
       &\cong&\projlim{ n} X_{nr}(k_{n,\infty}).
 \end{eqnarray*}

Here $Cl(K)(p)$ denotes the $p$-part of the ideal class group of a
(possibly infinite) extension $K$ of $\mathbb{Q}.$ Y. Ochi has
proven that $X_{nr}$ is a finitely generated torsion
$\Lambda(G)$-module.

\begin{thm}\label{xnr} Let $\mu_p\subseteq k$ and $k_\infty=k(A_{p^\infty})$ as above and suppose
that there is only one prime $\p$ of $k$ above $p$ and that $\p$
is totally ramified in $k_\infty/k.$ If $k_\infty/k$ is unramified
outside $\p,$ then there are canonical isomorphisms
\[(X_{nr})_{H_n}\cong X_{nr}(k_{n,\infty})\] and
\[(X_{nr})_{G_n}\cong Cl(k_{n,n})(p).\]

If, in addition, the $\mu$-invariant $\mu(X_{nr}(k_{cycl}))=0$ is
zero, then the following holds:
\begin{enumerate}
\item The $\mu$-invariant is zero for any intermediate field $\mu(X_{nr}(k_{n,\infty}))=0$
and $X_{nr}$ is a finitely generated $\Lambda(H)$-module. In
particular, the image of $X_{nr}$ in the quotient category is a
(possibly zero) completely faithful, hence cyclic object.
\item  There is an injective homomorphism of $\Lambda(G)$-modules \[Cl(k_\infty)(p)^\vee\to
\E^1(X_{nr})\] with pseudo-null cokernel, i.e.\ the Pontryagin
dual of the direct limit of the $p$-ideal class groups is
pseudo-isomorphic to the Iwasawa-adjoint of the inverse limit. In
particular, if $X_{nr}$ is pseudo-null, then   the ideal class
group $Cl(k_\infty)(p)=0$ vanishes.
\item  There exists a constant $c\geq 0$ independent of $n$ such that
\[\lambda(X_{nr}(k_{n,\infty}))=\rk_{\Lambda(H)}(X_{nr}) p^n +c\] holds for all sufficiently large
$n.$
\end{enumerate}
\end{thm}

We should mention that until now no single example of $X_{nr}$ not
being pseudo-null is known. On the contrary, at least  for the
first irregular prime $p=37,$ W.G. McCallum and R.T. Sharifi show
in \cite{MS} that for $k=\Q(\mu_{37})$ and
$k_\infty=k(\mu_{37^\infty}, 37^{37^{-\infty}})$ the Greenberg
conjecture holds, i.e.\ $X_{nr}(k_\infty)$ is a pseudo-null
$\La(G)$-module. Presumably, this is true for all primes, because
otherwise the set of irregular primes would be divided into two
classes corresponding to whether the Greenberg conjecture for the
above extension holds or not.

Before we can prove the theorem we need some preparation. For a
discrete $G_k$-module $N$ we set $N^\ast:=\Hom_\zp(N,\qp/\zp)$
where this module is endowed with the discrete topology.

\begin{lem}\label{e1}
For every $\Lambda$-module $M$ there is a natural exact sequence
\[0\to (\dirlim{n} M_{G_n}(p))^\vee\to \E^1(M)\to
(\dirlim{n}\H^1(G_n,M^\vee)^\ast\otimes\qp/\zp)^\vee\to 0.\]

In particular, if $M$ is pseudo-null, then the direct limit
$\dirlim{n} M_{G_n}(p)=0$ vanishes.
\end{lem}

\begin{proof}
This is just the Pontryagin dual of \cite[(5.4.13)(i)]{nsw} for
$r=1.$
\end{proof}

\begin{prop}
Let $M$ be a $\Lambda(G)$-module which is finitely generated over
$\Lambda(H).$
\begin{enumerate}
\item Assume that $M_{G_n}$ is finite.  Then $(M_{H_n})^{\Gamma_n}$ is finite, too.
\item Assume that $M$ has no non-zero pseudo-null submodule. Then it holds that
\[\H^1(G_n,M^\vee)\cong ((M_{H_n})^{\Gamma_n})^\vee.\]
If, in addition, $M_{G_n}$ is finite for all $n,$ then we obtain
an isomorphism
\[(\dirlim{n} M_{G_n})^\vee\cong \E^1(M).\]
\item Assume again that $M_{G_n}$ is finite for all $n.$ Then there is a canonical injection
\[(\dirlim{n} M_{G_n})^\vee\to E^1(M)\] with pseudo-null cokernel.
\end{enumerate}
\end{prop}

\begin{proof}
The first item follows from the exact sequence
\[0\to (M_{H_n})^{\Gamma_n} \to M_{H_n} \to   M_{H_n} \to M_{G_n} \to 0\]
noting that $M_{H_n}$ is a finitely generated \zp-module. Now
assume that $M$ has no non-zero pseudo-null submodule, i.e.\ $M$
is a torsion-free $\Lambda(H)$-module (cf. \cite[lem 6.4
]{ochi-ven}).  Thus $(M^{H_n})_{\Gamma_n}=0$ and the first part of
(ii) follows from the Hochschild-Serre long exact sequence

\[ 0\to ((M_{H_n})^{\Gamma_n})^\vee\to \H^1(G_n, M^\vee)\to ((M^{H_n})_{\Gamma_n})^\vee=0\]
while the second part is a consequence of the first part, (i) and
lemma \ref{e1}. To prove (iii) we denote by $M_0$ the maximal
pseudo-null submodule of $M$ and we set $N:=M/M_0.$ The exact
sequence \[0\to M_0\to M \to N\to 0\] induces via   the long exact
sequence of $\E^i$ the exact sequence \begin{eqnarray}
\label{e1-iso} \E^0(M_0)=0\to \E^1(N)\to \E^1(M)\to
0=E^1(M_0)\end{eqnarray} using the fact that
$\E^0(M_0)=\E^1(M_0)=0$ vanishes for pseudo-null modules. On the
other hand we obtain taking co-invariants \[0=N^{H_n}\to
(M_0)_{H_n} \to M_{H_n} \to N_{H_n} \to 0\] and \[
(N_{H_n})^{\Gamma_n}\to (M_0)_{G_n} \to  (M)_{G_n} \to (N)_{G_n}
\to  0.\] First note that $ (N)_{G_n}$ and hence
$(N_{H_n})^{\Gamma_n}$ are finite for all $n$ by (i). Thus also
$(M_0)_{G_n}$ is finite for all $n$. So, by lemma \ref{e1},
$\dirlim{n} M_{G_n}\cong\dirlim{n} N_{G_n}$ and the result follows
from the above isomorphism \ref{e1-iso} and (ii).
\end{proof}

\begin{proof} (of the theorem)
The first statement is proved by the same arguments as for
$\mathbb{Z}_p$-extensions, see e.g.\ \cite[lem. 11.1.5]{nsw}. Now
assume that the $\mu$-invariant vanishes.  Nakayama's lemma tells
that $X_{nr}$ is a finitely generated $\Lambda(H)$-module, i.e.\
$X(k_{n,\infty})=(X_{nr})_{H_n}$ a finitely generated $\zp$-module
and thus $\mu(X(k_{n,\infty}))=0.$ Item (ii) follows from the
above proposition while (iv) is a consequence of the next lemma.
\end{proof}

\begin{lem}
Let $M$ be a $\Lambda(G)$-module which is finitely generated as
$\Lambda(H)$-module of rank $\rk_{\Lambda(H)}M=d.$ Then there is a
constant $c\geq 0$ such that \[\rk_\zp M_{H_n}=dp^n +c\] for all
sufficiently large $n.$
\end{lem}

\begin{proof}
Since $M$ is finitely generated as
$\Lambda:=\Lambda(H)\cong\zp[[X]]$-module, by the (commutative)
structure theory we obtain an exact sequence
\[0\to M/M_0\to \Lambda(H)^d\oplus\bigoplus_{i=0}^s \Lambda/\Lambda F_i\oplus\bigoplus_{j=0}^r\Lambda/\Lambda p^{n_j}\to E\to 0,\]
where $M_0$ denotes the maximal finite submodule of $M,$ $E$ is
finite, the $n_j\geq 1,$ $0\leq r,$ are natural numbers and the
$F_i\in\Lambda,$ $0\leq i\leq s,$ are Weierstrass polynomials.
Taking $H_n$-coinvariants and tensoring with $\qp$ one sees
immediately that \[\rk_\zp
M_{H_n}=\rk_\zp\zp[H/H_n]^d+\sum_{i=0}^s \rk_\zp (\Lambda /\Lambda
F_i)_{H_n}\] holds. Since the second summand is bounded by
$\sum_{i=0}^s\rk_\zp \Lambda /\Lambda F_i<\infty,$ it becomes
stationary  for large $n$ and the result follows.
\end{proof}

The asymptotic behaviour of the $\lambda$-invariant suggest the
following  asymptotic growth of the order of the $p$-ideal class
group. If $e_n$ denotes the precise exponent of
$\#Cl(k_{n,n})(p)=p^{e_n}$ with respect to $p,$ then I guess that

\[e_n=(\rk_{\Lambda(H)}X_{nr} \cdot n + O(1))p^n\] holds.

We finish this paper with a corollary on certain units which are
linked to the module $X_{nr}$ via Kummer theory. To this end, we
denote by $\mathcal{O}_K$ the ring of integers of a number field
$K$ and its units by $E(K):=\mathcal{O}_K^\times.$ We set
$E(k_\infty)=\dirlim{} E(k')$ where the direct limit runs through
all finite subextensions $k|k'|k_\infty.$ Then
\[\mathcal{E}(k_\infty):=(E(k_\infty)\otimes_\mathbb{Z}\qp/\zp)^\vee\] is a compact
 $\La(G)$-module which fits into the following exact and commutative diagram (see \cite[prop.
 4.15]{ven})

 \begin{eqnarray*}
\xymatrix{
    {\ 0 \ar[r] } &  {\ Cl(k_\infty)(p)^\vee \ar[r] } &  {\  X_\Sigma(-1) \ar[r] } &  {\ \e(k_\infty) \ar[r] } &  {\
    0, } \\
    {\ 0 \ar[r] } &  {\ Cl(k_\infty)(p)^\vee \ar[r]\ar@{=}[u] } &  {\ \tor_\La X_\Sigma(-1) \ar[r]\ar@{^(->}[u] } &  {\ \tor_\La\e(k_\infty) \ar[r]\ar@{^(->}[u] } &  {\
    0. }
 }&&
\end{eqnarray*}

Here, $\Sigma:=S_p\cup S_\infty$ and
$X_\Sigma(-1)=\H^1(G_\Sigma(k_\infty),\mu_{p^\infty})^\vee$ is the
first negative Tate twist of the Galois group
$G_\Sigma(k_\infty)^{ab}(p)$   of the  maximal outside $\Sigma$
unramified abelian $p$-extension $k_\Sigma^{ab}(p)$ over
$k_\infty.$ With other words, $G_\Sigma(k_\infty)$ is the Galois
group of the maximal outside $\Sigma$ unramified extension
$k_\Sigma$ of $k$ over $k_\infty$ and
$G_\Sigma(k_\infty)^{ab}(p)$ its maximal abelian pro-$p$ quotient.
It follows from \cite[thm 3.1.5, cor. 3.1.6]{ven}, that $\tor_\La
X_\Sigma(-1)$ is pseudo-isomorphic (in the quotient category) to
$\E^1(X_{nr}).$ Finally, as a consequence of the above theorem we
obtain

\begin{cor}
Under the assumptions of the theorem the torsion submodule
$\tor_\La\e(k_\infty)$ of $ \e(k_\infty)$ is pseudo-null.
\end{cor}

\providecommand{\bysame}{\leavevmode\hbox
to3em{\hrulefill}\thinspace}
\providecommand{\MR}{\relax\ifhmode\unskip\space\fi MR }
\providecommand{\MRhref}[2]{%
  \href{http://www.ams.org/mathscinet-getitem?mr=#1}{#2}
} \providecommand{\href}[2]{#2}

\end{document}